\documentclass[12pt]{article}
\usepackage{amsfonts,amssymb}
\usepackage{graphicx} 
\usepackage{float}
\usepackage{color}
\usepackage{amsthm, amssymb, latexsym, amsmath, bm %showkeys
}
\usepackage[margin=1.35 in]{geometry}
\usepackage{enumerate}

\usepackage{dsfont}

\usepackage{hyperref}
\usepackage{cleveref}
\usepackage[dvipsnames]{xcolor}

\usepackage[font=small,labelfont=bf]{caption}
\usepackage{subcaption}

\hypersetup{
	colorlinks   = true, %Colours links instead of ugly boxes
	urlcolor     = blue, %Colour for external hyperlinks
	linkcolor    = blue, %Colour of internal links
	citecolor   = red %Colour of citations
}

\newtheorem{theorem}{Theorem}[section]

\newtheorem{lemma}[theorem]{Lemma}
\newtheorem{proposition}[theorem]{Proposition}

\theoremstyle{definition}

\numberwithin{theorem}{section}
\numberwithin{equation}{section}

\newcommand{\Del}{\Delta}
\newcommand{\eps}{\varepsilon}

\newcommand{\R}{{\mathbb R}}

\newcommand{\bC}{{\mathbb C}}

\newcommand{\cF}{\mathcal{F}}

\newcommand{\cM}{\mathcal{M}}
\newcommand{\cN}{\mathcal{N}}
\newcommand{\cS}{\mathcal{S}}
\newcommand{\cU}{\mathcal{U}}
\newcommand{\cW}{\mathcal{W}}

\newcommand{\Ran}[1]{\operatorname{Ran}\, #1}

\renewcommand{\Re}{\operatorname{Re}}
\renewcommand{\Im}{\operatorname{Im}}

\newcommand{\dist}{{\rm{dist\, }}}

\newcommand{\p}{\partial}

\newcommand{\DETAILS}[1]{}

\newcommand{\variants}[1]{}

\newcommand{\N}{\mathbb{N}_0}
\begin{document}
	\title{Stability of Traveling Fronts of the FitzHugh–Nagumo Equations on Cylindrical Surfaces}

	\author{Afroditi Talidou*}
    \date{\small{*Hotchkiss Brain Institute, University of Calgary, Calgary, AB T2N 4N1, Canada}} 

	\maketitle
	%\vspace{-2 cm}

\begin{abstract} 

The FitzHugh-Nagumo equations are known to admit traveling front solutions in one spatial dimension that are nonlinearly stable. This paper concerns the stability of traveling front solutions propagating on cylindrical surfaces. It is shown that such traveling fronts are nonlinearly stable on the surface of standard cylinders of constant radius. The analysis is extended to warped cylinders with slowly varying radius, where persistence of front-like solutions is established. Numerical simulations support the theoretical findings.

\end{abstract}

\section{Introduction}\label{sec:introduction} 

Neuronal axons are elongated projections of neurons responsible for transmitting electrical signals, known as action potentials, from the cell body to other neurons or muscles. Their geometry -- characterized by long, thin cylindrical structures -- is essential for the efficient conduction of these signals over extended distances. Mathematical models of action potential propagation in axons commonly employ reaction-diffusion equations, which capture the dynamics of ion movement and membrane potential changes underlying the generation and transmission of action potentials.

Characterized as ultrathin cables, neuronal axons are typically modeled as straight lines without internal structure. In this work, we instead represent axons as thin cylindrical surfaces. Let $\cS$ denote such a surface. On $\cS$, we consider the membrane potential $u = u(x, \theta, t)$ as a function of the longitudinal spatial variable $x$, the angular coordinate $\theta$, and time $t$. The cable equation governing the membrane potential is
\begin{equation}\label{eq:cable}
C_\mathrm{m} \p_t u = G\,\Delta_{\cS}u - I_m\,,
\end{equation}   
where $C_\mathrm{m}$ is the membrane capacitance per unit length, $G$ is the conductance and $I_m$ is the ionic membrane current density specified by the chosen reaction kinetics \cite[Section 6.3]{Abbott_book}. The conductance per unit length is equal to the inverse of the resistance and is given by
\begin{equation}\label{eq:G}
G := \frac{\pi R^2}{r_{\mathrm{int}}}\,.
\end{equation}
Here, $R$ is the radius of the axon and $r_{\mathrm{int}}$ is the intracellular resistivity.

For $I_m$ we choose the FitzHugh-Nagumo (FHN) dynamics \cite{FH61, Nagumo62}, a simplification of the Hodgkin-Huxley dynamics \cite{HH52}, which describe the neuronal excitability and recovery of the membrane potential. As such, Eq.~\eqref{eq:cable} takes the form
\begin{equation}\label{eq:FHN}
\begin{split}
C_\mathrm{m} \p_t u_1 &= G\,\Delta_{\cS}u_1 + f(u_1) - u_2\,, \\
\p_t u_2 &= \eps (u_1 - \gamma u_2)\,,
\end{split}
\end{equation}
where $\Delta_{\cS}$ is the Laplace-Beltrami operator on a cylindrical surface $\cS$. In Eq.~\eqref{eq:FHN}, $u_1$ represents the excitable behaviour of the membrane potential, that is, the rapid response to a stimulus that exceeds a certain threshold, and $u_2$ is a recovery variable (acting at a slower time scale) that allows the membrane potential to return to its resting state and become excitable again. The nonlinear function $f$ is the positive feedback and is defined by a cubic polynomial $f(u_1) = -u_1(u_1 - \alpha)(u_1 - 1)$, where $\alpha \in \left( 0, \frac{1}{2} \right)$ is the threshold of neuronal firing. The parameters $\eps$ and $\gamma$ determine the timescales of the system: $\eps$ controls the rate at which the membrane potential $u$ returns to its resting state (typically assumed to be small ($\eps \ll 1$)), while $\gamma$ controls the responsiveness of the recovery variable $w$ to changes in membrane potential.

Equation~\eqref{eq:FHN} has a traveling wave solution of the form $(u_1, u_2) = (\phi_1(z), \phi_2(z))$, $z=x-ct$, $c>0$ that satisfies
\begin{equation}\label{eq:front}
\begin{split}
\lim_{z\to -\infty} (\phi_1(z), \phi_2(z)) &= (1, 1)\,,\\
\lim_{z\to +\infty} (\phi_1(z), \phi_2(z)) &= (0, 0)\,.
\end{split}
\end{equation}
This is called a traveling front solution, or simply a front. 

Taking $\cS = \mathbb{R}$, the system \eqref{eq:FHN} describes solutions on the real line. For this special case, early rigorous results established the existence of traveling waves (including traveling fronts) and their stability. The existence of traveling front solutions of Eq.~\eqref{eq:FHN} on $\cS = \R$ has been proven both analytically, using singular perturbation theory, and numerically \cite{Rinzel_Terman}. It was shown in \cite{Bell_Deng} that the FHN equation admits infinitely many distinct traveling front and back solutions. The stability of the fronts was established by Jones~\cite{Jones_84} through singular perturbation and Evans function techniques, and independently by Yanagida~\cite{Yanagida_fronts} via eigenvalue analysis. Bifurcation theory and topological methods have also been employed to study the existence and stability of $N$-front solutions \cite{Nii, Sandstede}. More recently, invasion fronts have been linked to wave-train selection phenomena \cite{Carter_Scheel}, and the nonlinear stability of invasion-induced periodic patterns has been demonstrated \cite{Avery_Carter_et_al}.

In a previous study, we proved that nearly planar traveling pulse solutions of Eq.~\eqref{eq:FHN}, perturbed by small fluctuations, propagate on a warped cylindrical surface $\cS$ \cite{Talidou_21}. On undulated cylinders, such solutions have been shown to behave comparably to their radial averages under specific conditions \cite{Karali}. Exact solutions of a generalized cable equation on tapered cylindrical geometries have also been studied \cite{Timofeeva_20}.

In this work, we adopt a framework similar to that of \cite{Talidou_21}, but extend it in two key directions. First, while the analysis in \cite{Talidou_21} focuses on pulse solutions of the FHN equations \eqref{eq:FHN} on cylindrical surfaces $\cS$, in the present study we investigate traveling front solutions of Eq.~\eqref{eq:FHN} on $\cS$. Second, in \cite{Talidou_21} the conductance $G$ and membrane capacitance $C_m$ were normalized to unity, whereas here we retain the full definition of the conductance (see Eq.~\eqref{eq:G}), which depends explicitly on the cylinder radius $R$. This $R$-dependence of $G$ is crucial for the stability analysis, as it appears in the leading term of the first equation in \eqref{eq:FHN}.

Within this framework, we investigate solutions of Eq.~\eqref{eq:FHN} on thin, infinitely long cylindrical surfaces of two types: (i) standard cylinders with constant radius, and (ii) warped cylinders whose radii vary slowly along their length, representing more realistic neuronal axon morphologies. Specifically, we show that traveling fronts exhibit the following behavior: on standard cylinders of small constant radius, they are stable under general perturbations of the initial condition, whereas on warped cylinders, solutions that are initially close to a traveling front remain close to the family of propagating fronts for all time.

We support our theoretical findings with numerical simulations. We provide an example of a cylinder of constant radius, followed by two examples of warped cylinders. The first one considers a cylindrical surface that closely resembles neuronal axons exhibiting the characteristic pearls-on-a-string axonal morphology \cite{pearl}. These axons consist of a series of tiny beads (“pearls”) connected by thin cylindrical segments (“string”). The second one considers a cylindrical geometry with a single swelling \cite{Kolaric}. For all three configurations, we identified parameter values for which a traveling front is generated and is stable as it propagates along the axon.

\paragraph{Main results.} We begin with the standard cylinder of constant radius $R$
\begin{equation*} \label{eq:def-SR}
\cS_R := \left \{ (x,R\cos\theta, R\sin\theta) 
\in\R^3\ \big\vert\ 
x\in \R,\, \theta\in [0, 2\pi)
\right\}.
\end{equation*}
On $\cS_R$, the Riemannian area element is $R d\theta dx$ and the Laplace operator takes the form 
\begin{equation}\label{eq:Delta_S_R}
\Delta_{\cS_R} = \p_x^2 + R^{-2}\p_\theta^2\,.
\end{equation}
Equation~\eqref{eq:FHN} is invariant under translations, that is, if $u(x, \theta, t)$ solves Eq.~\eqref{eq:FHN}, then so do its translates $u_h(x, \theta, t) := u(x - h, \theta, t)\,, h\in \R$. A front $\Phi$, traveling with constant speed $c$, is an axisymmetric traveling wave solution $u(x,\theta,t) = \Phi(x-ct)$ of Eq.~\eqref{eq:FHN} on $\cS_R$. By translation invariance, the entire family of translates $\Phi(x-h-ct)$, $h\in \mathbb{R}$, are also solutions of the same speed $c$.

Our analysis begins with the stability of such fronts under perturbations of the initial values. To formulate this precisely, we work in the mixed Sobolev space
\begin{equation}\label{eq:H21}
H^{2,1}_R
:=\left\{u=(u_1, u_2)\in L^2\,\times\,L^2\ \big\vert\ 
\Delta_{\cS_R} u_1\in L^2, \p_x u_2\in L^2\right\}\,,
\end{equation}
equipped with the norm $\| u \|_{2,1} := \| (\Delta_{\cS_R} u_1, \p_x u_2) \| + \| u \|$. Here $\| \cdot \|$ denotes the $L^2-$norm. With this setting, the initial value problem for Eq.~\eqref{eq:FHN} is locally well-posed. The proof is similar to \cite[Proposition 2.1]{Talidou_21} and is omitted.

The translates of $\Phi$ form a one-dimensional manifold of front solutions
\begin{equation}\label{eq:manifold}
\cM := \{ \Phi_h \;\vert\; h\in \R \}\,.
\end{equation}
The following theorem indicates that $\cM$ is asymptotically orbitally stable in $\| \cdot \|_{2,1}$ with exponential rate $\mu >0$.

\begin{theorem}[Stability of fronts on $\cS_R$]
	\label{thm:nonli-stab}
    Let $\cS_R$ be a standard cylinder of constant radius $R\leq 1$ and consider Eq.~\eqref{eq:FHN} on that geometry. Assume $0< \eps \ll 1$, $\gamma>0$ and $\alpha \in \left( 0, \frac12\right)$ are such that Eq.~\eqref{eq:FHN} has a front solution $\Phi(x-ct)$, and that $C_\mathrm{m}, r_{int} >0$. Then, there exists a neighborhood $\cU$ of the front solution $\Phi(x)$ in $H^{2,1}_R$ such that for every initial value $u_0\in \cU$, the mild solution $u(t)$ of Eq.~\eqref{eq:FHN} with initial condition $u\big\vert_{t=0}=u_0$
	exists globally in time. Moreover, it satisfies
	\begin{equation}\label{eq:nonli-stab}
	\|u(t)- \Phi_{ct+\tilde{h}} \|_{2,1} 
	\le K_1 e^{-\mu t} \|u_0-\Phi\|_{2,1}\,, \qquad (t\ge 0)
	\end{equation}
	for $K_1 > 0$ independent of $\eps$,  $\mu \gtrsim\eps$, and $\tilde{h}\in\R$ 
	with $$|\tilde{h}|\le K_2\|u_0-\Phi\|_{2,1}^2\,,
	$$
    where $K_2>0$ and independent of $\eps$.
\end{theorem}
\noindent

The proof of Theorem~\ref{thm:nonli-stab} relies on \cite[Theorem 4.3.5]{KPbook} and is given in Section \ref{sec:front_standard_cylinder}. For convenience, we have included \cite[Theorem 4.3.5]{KPbook} in Appendix~\ref{app:stability_thm}. The analysis combines functional analytic techniques and spectral theory to derive resolvent estimates extending the frameworks developed in \cite{Talidou_21}.

The distance of $u\in H_R^{2,1}$ from the manifold $\cM$ is $\dist(u, \cM) := \inf_h \| u - \Phi_h \|_{2,1}$. Theorem~\ref{thm:nonli-stab} implies that there is a tubular neighbourhood
$\cW = \{w\in H_R^{2,1} \mid \dist(w,\cM)<\xi\}$ of $\cM$ for which
\begin{equation*}
\dist(u(t), \cM)\le C_1e^{-\mu t} \dist(u_0, \cM)
\end{equation*}
for all mild solutions with initial values in $\cW$. As $t\to\infty$, each solution converges to a particular front $\Phi(x-ct-\tilde{h})$. This holds by translation invariance.

Turning to warped cylindrical surfaces $\cS_{\rho}$, we define them as graphs over the standard ones determined by a positive function $\rho(x)$:
\begin{equation}\label{eq:def-Srho}
\cS_\rho := \left \{ x\in\R,\; \theta\in[0, 2\pi) \;\;\big\vert\;\; (x, \rho(x)\cos\theta, \rho(x)\sin\theta) \in \R^3 \right\}\,.
\end{equation}
The Laplace-Beltrami operator on $\cS_\rho$ is given by
\begin{equation}\label{eq:Delta-rho}
\Del_{\cS_{\rho}} :=
\frac{\pi}{r_{\mathrm{int}}}\left[ \frac{1}{\sqrt{g(x)}} \p_x \left( \frac{\rho^3(x)}{\sqrt{1+\rho'(x)^2}} \,\p_x \right) + \p_{\theta}^2 \right]\,,
\end{equation}
where $g(x) := \rho(x)^2(1 + \rho'(x)^2)$.

Assuming that $\rho \in C^2$ is positive, bounded, and bounded away from zero, we define the diffeomorphism $\psi_{\rho} \in C^2$ such that
\begin{equation}\label{eq:diffeo}
    \psi_\rho(x,R\cos\theta, R\sin\theta)
    = \bigl(x, \rho(x)\cos\theta,\rho(x)\sin\theta\bigr)\,.
\end{equation}
Functions on $\cS_{\rho}$ are identified as functions on $\cS_R$ via $\psi_{\rho}$, hence the norms $\| \cdot \|$ and $\| \cdot \|_{2,1}$ on the standard cylinder of constant radius are applied to functions on $\cS_{\rho}$.

As in the case of traveling pulses studied in \cite{Talidou_21}, when the radius $\rho$ is non-constant, Eq.~\eqref{eq:FHN} does not possess traveling front solutions. Our next result shows, however, that front-like solutions persist when $\rho$ is almost constant. Moreover, these front-like solutions remain close to the manifold $\cM$ and evolve along $\cM$ over time.

\begin{theorem}[Persistance of front-like solutions on $\cS_{\rho}$]\label{thm:persist}
    Let $\cS_{\rho}$ denote the warped cylinder of variable radius defined in Eq.~\eqref{eq:def-Srho}. Consider Eq.~\eqref{eq:FHN} on $\cS_{\rho}$ with parameters $\alpha$, $\eps$, $\gamma$, $C_\mathrm{m}$, and $r_{\mathrm{int}}$ as in Theorem~\ref{thm:nonli-stab}. Under these assumptions, there exists a positive constant $\delta_*$ and a tubular
	neighborhood $\cW$ of $\cM$ in $H_R^{2,1}$ such that if $\delta:=R^{-1}\|\rho - R\|_{C^2}\le \delta_*$, with $R \leq 1$, then for every initial value $u_0\in\cW$, the unique mild
	solution $u(t)$ with $u\big\vert_{t=0}=u_0$
	exists globally in time, and satisfies
	\begin{equation}\label{eq:persist}
	\dist(u(t), \cM) \leq K'_1 e^{-\nu t}\dist(u_0,\cM) + K'_2\delta
	\qquad (t\ge 0)\,.
	\end{equation}
	The above constants are given explicitly by: $K'_1 = 2K_1$, $K'_2 = (2+K_1)C_0$ and $\nu = \frac{\mu \ln2}{ln(2K_1)}$. The constants $K_1$ and $\mu$ are defined in Theorem~\ref{thm:nonli-stab}, $C_0$ is as in Proposition~\ref{prop:continuous-rho} and $T=\mu^{-1} \log(2K_1)$.
\end{theorem}

Under the assumptions of the theorem, initial values $u_0\in\cM$ give rise to front-like solutions that satisfy 
$\sup_{t>0} \dist(u(t), \cM)\le K'_2\delta$. The proof of Theorem~\ref{thm:persist} is given in Section \ref{sec:warped}. The tubular neighborhood $\cW$, as well as the constants $\delta_*$, $K'_1$, $K'_2$, and $\nu$ depend on the parameters $\alpha$, $\gamma$, $C_\mathrm{m}$ and $r_{\mathrm{int}}$.

\paragraph{Outline.} In Section~\ref{sec:front_standard_cylinder} we prove Theorem~\ref{thm:nonli-stab}. We begin by establishing the linearized stability. Consider Eq.~\eqref{eq:FHN} on the standard cylinder of constant radius $R$. In the moving frame $z=x-ct$, the equation takes the form
\begin{equation*}
    \partial_t u = F(u),
\end{equation*}
where $F(u)$ denotes the right-hand side of Eq.~\eqref{eq:FHN_moving} given below. A Taylor expansion of $F$ about $\Phi$ splits it into its linear and nonlinear components. Denote by $L$ the linearization about the stationary solution $\Phi$. The main difficulty lies in showing that the semigroup generated by $L$ decays exponentially in directions transversal to the tangent space of $\cM$ at $\Phi$. Specifically, we seek to prove the estimate
\begin{equation*}
    \left\| e^{tL} (I-P) \right\|_{2,1} \lesssim e^{-\eta t} \;\;\; (t\geq 0)\,,
\end{equation*}
for $\eta := \min\left\{ C_m^{-1} |f'(\phi_1^{\pm})|, \beta, \eps \gamma \right\}$, where $\beta$ is the spectral gap given in Lemma~\ref{lem:L0}. The operator $P$ is a projection onto the tangent space that commutes with $L$, and $\| \cdot \|_{2,1}$ is the operator norm on $H^{2,1}$. Since $L$ is not self-adjoint, $P$ is not an orthogonal projection. The tangent space of $\cM$ at $\Phi$ is spanned by the derivative $\tau:=-\partial_z\Phi$, and $L\tau = 0$. Hence, $P$ is constructed as the spectral projection associated with the zero eigenvalue of $L$. The decay estimate follows since the remainder of the spectrum lies in the left half-plane. Nonlinear stability can then be deduced by direct application of Theorem 4.3.5 in \cite{KPbook}.

In Section~\ref{sec:warped} we prove Theorem~\ref{thm:persist}. Linearizing about the traveling front $\Phi$ in the moving frame, as in the proof of Theorem~\ref{thm:nonli-stab}, yields a time-dependent perturbation of the principal part of the linearized operator, because $\Delta_{\cS_{\rho(x)}}$ becomes $\Delta_{\cS_{\rho(z+ct)}}$. The resulting operator is not self-adjoint, and estimating the evolution it generates is delicate. To address this difficulty, we consider instead the flow on a warped cylinder $\cS_{\rho}$ in the static frame as a perturbation of the flow on the standard cylinder $\cS_R$. We estimate the difference by first comparing the flows linearized around $u=0$, and then applying Gronwall's inequality to extend this estimate to the nonlinear flow. This approach is relatively simple, since the linearized operator generates an analytic semigroup on $L^2$ that restricts to a uniformly bounded analytic semigroup on the dense invariant subspace $H_R^{2,1}$. Combining this estimate for the flows on $\cS_{\rho}$ and $\cS_R$ with the exponential decay of fluctuations for front solutions on $\cS_R$, established in Theorem~\ref{thm:nonli-stab}, yields Theorem~\ref{thm:persist}. 

In Section~\ref{sec:numerics} we present numerical simulations for both the standard and warped cylinders, and in Section~\ref{sec:discussion} we conclude with a discussion of future work.

%%%%%%%%%%%%%%%%%%%%%%
\section{Front solutions on standard cylinders}
\label{sec:front_standard_cylinder}

We consider first the case of a standard cylinder $\cS = \cS_R$ of constant radius $R$. Fix a front solution $\Phi$ on the surface of $\cS_R$, and let $\mathcal{M}$ be the manifold of its translates defined in Eq.~\eqref{eq:manifold}. By definition, $\Phi$ is an axisymmetric traveling wave solution of the FHN system \eqref{eq:FHN}. Adopting a coordinate system defined by \( z = x - ct \), the system in Eq.~\eqref{eq:FHN} transforms into the following 
\begin{equation} \label{eq:FHN_moving}
\begin{aligned}
\partial_t u_1 &= D \Delta_{\mathcal{S}_R} u_1 + c \partial_z u_1 + C_\mathrm{m}^{-1}(f(u_1) - u_2), \\
\partial_t u_2 &= c \partial_z u_2 + \eps (u_1 - \gamma u_2),
\end{aligned}
\end{equation}
where \( \Phi \) becomes a stationary solution in this frame. The effective diffusivity \( D \) is
\[
D := C_\mathrm{m}^{-1}\,G  = C_\mathrm{m}^{-1}\frac{\pi R^2}{r_{\mathrm{int}}}, \quad \text{with } C_\mathrm{m}, r_{\mathrm{int}} > 0.
\]
Let $F(u)$ be the right hand side of Eq.~\eqref{eq:FHN_moving}. Since \( \Phi \) is a stationary solution, it satisfies \( F(\Phi) = 0 \). Set $u := \Phi + v$, where $\Phi\in \cM$ and $v$ a fluctuation. A Taylor expansion about $\Phi$ yields
\begin{equation}
\p_t v = F(\Phi + v) = Lv + N(v)\,,
\end{equation}
where the linearization $L$ is given by the G\^ateaux derivative of $F$ about $\Phi$:
\begin{equation} \label{eq:L}
L := dF(\Phi) =
\begin{pmatrix}
D \Delta_{\mathcal{S}_R} + c \partial_z + C_{\mathrm{m}}^{-1}f'(\phi_1) & -C_{\mathrm{m}}^{-1} \\
\eps & c \partial_z - \eps \gamma
\end{pmatrix},
\end{equation}
and the nonlinearity $N(v)$ is defined as $N(v) := F(\Phi + v) - Lv$, or explicitly 
\begin{equation}\label{eq:N(v)}
N(v) :=
\begin{pmatrix}
C_{\mathrm{m}}^{-1}\,v_1^2\,(-v_1 - 3\phi_1 + \alpha + 1)\\
0
\end{pmatrix}\,.
\end{equation}

The linearization in Eq.~\eqref{eq:L} defines a closed, linear operator on the Hilbert space $L^2 := L^2(\cS_R; \mathbb{C}^2)$ of two-component square integrable functions. The inner product of this space is weighted to reflect the separation of time scales:
\begin{equation} \label{eq:eps-inner}
\langle u, w \rangle := \int_{\mathbb{R}} \int_{S^1} \left( u_1 \bar{w}_1 + \eps^{-1} u_2 \bar{w}_2 \right) R\, d\theta\, dz.
\end{equation}
The domain of $L$ is $H^{2,1}_R$ defined in Eq.~\eqref{eq:H21} with the norm $\| \cdot \|_{2,1}$. This norm is equivalent to the graph norm of $L$ (cf. Lemma~\ref{lem:L-domain}). As the angular derivative $\p_{\theta} u_2$ is not required to be in $L^2$, the space $H_R^{2,1}$ properly contains \( H^2 \times H^1 \). We will show below that \( L \) generates a strongly continuous semigroup \( e^{tL} \) on this domain, which is essential for the subsequent stability analysis.

Due to translational invariance of the governing equations, the family of shifted profiles \( \Phi_h \) also satisfies \( F(\Phi_h) = 0 \) for all \( h \in \mathbb{R} \). Differentiating this identity with respect to \( h \) yields the relation $L \p_z \Phi = 0$, indicating that 
\begin{equation}\label{eq:span}
\operatorname{Ker}(L) = \mathrm{span}\{ \p_z \Phi \}\,,
\end{equation}
and that zero is an eigenvalue of the linearized operator. 
%neutral direction in dynamics

%%%%%%%%%%
%%%%%%%%%%
%%%%%%%%%%

\subsection{The semigroup generated by $L$}
\label{subsec:semigroup}

Consider the linearization $L$, about the front solution $\Phi$, defined in Eq.~\eqref{eq:L} with domain $H_R^{2,1}$. The main result of this subsection is stated next.

\begin{proposition}\label{prop:lin_semigroup}
	The linearized operator $L$ generates a strongly continuous semigroup, denoted as $e^{tL}$, satisfying the estimate
	\begin{equation}\label{eq:semigroupL}
		\| e^{tL} \| \leq e^{t}\,, \qquad \text{for all } t \geq 0\,.
	\end{equation}
\end{proposition}

The idea is to write $L$ as a generator of a strongly continuous semigroup, that is easy to study, perturbed by a bounded operator. To that end, decompose $L$ as $L = L_{\pm} + V_{\pm}$. The operators $L_{\pm}$, defined as 
\begin{equation}\label{eq:Lpm}
L_{\pm} :=
\begin{pmatrix}
D \Delta_{\cS_R} + c\p_z + C_{\mathrm{m}}^{-1} f'(\phi_1^{\pm}) & -1 \\
\eps & c\p_z - \eps \gamma
\end{pmatrix}\,,
\end{equation} 
describe the asymptotic behavior of $L$ as $z\to \pm \infty$. Specifically, the wave profile $\Phi$ satisfies $\lim_{z\to +\infty} \Phi(z) = \phi^+_1 = 0$ and $\lim_{z\to -\infty} \Phi(z) = \phi^-_1 = 1$. The matrix multiplication operators, defined by
\begin{equation}
V_{\pm} := 
\begin{pmatrix}
C_\mathrm{m}^{-1}\left(f'(\phi_1(z)) - f'(\phi_1^\pm)\right) & 1 - C_\mathrm{m}^{-1} \\
0 & 0
\end{pmatrix}\,,
\end{equation}
are bounded on $L^2$ and $H^{2,1}$ because $f$ is a polynomial and $\phi_1$ is a smooth, bounded function. Moreover, the first entry of $V_{\pm}$ is continuous  and decays at infinity. By (Theorem 3.1.11, \cite{KPbook}) the operator $L$ is a bounded, relatively compact perturbation of $L_{\pm}$.

The following lemma implies that $L_{\pm}$ are dissipative operators in $L^2$, that is, $\mathrm{Re} \langle L_{\pm}v, v \rangle \leq 0$ on $H^{2,1}$, and as such they generate strongly continuous semigroups of contractions (cf. Lemma \ref{lem:dissipative}).

\begin{lemma}\label{lem:L_dissip}
	Let $L_{\pm}$ be the constant-coefficient operators defined in Eq.~\eqref{eq:Lpm}. Then, for all $v \in H_R^{2,1}$, $ \mathrm{Re} \langle L_{\pm} v, v \rangle \leq -\eta \|v\|^2 $ with $\eta := \min\left\{ C_\mathrm{m}^{-1}|f'(\phi_1^\pm)| ,\, \eps \gamma \right\}$.
\end{lemma}
\begin{proof}
	By the definition of the inner product \eqref{eq:eps-inner} we have
	\begin{align*} 
	\Re\, \langle L_{\pm} v, v \rangle
	&= \int_{\cS_R}
	\left((D \Delta_{\cS_R}v_1)\bar v_1 + C_\mathrm{m}^{-1} f'(\phi_1^{\pm}) |v_1|^2 - 
	\gamma |v_2|^2\right)\, R \,d\theta dz \notag\\
	&= -\int_{\cS_R}
	\left(C_\mathrm{m}^{-1} \left|f'(\phi_1^{\pm}) \right| |v_1|^2 + \gamma |v_2|^2\right)\, R\, d\theta dz\\
	&\le -\min\{C_\mathrm{m}^{-1} \left|f'(\phi_1^{\pm}) \right|, \eps\gamma\} \| v \|^2\,. 
	\end{align*}
\end{proof}

The operators $L_{\pm}$ possess only essential spectrum, as in the Fourier representation they are given by the matrix multiplication operators
\begin{equation}
\label{eq:multOp}
m_{\pm}(k,n) :=
\begin{pmatrix}
- D(k^2 + R^{-2} n^2) + i c k + C_\mathrm{m}^{-1} f'(\phi_1^\pm) & -1 \\
\eps & i c k - \eps \gamma
\end{pmatrix}\,.
\end{equation}
The resolvent set consists of all $\lambda \in \mathbb{C}$ for which the matrix \( \lambda - m_{\pm}(k, n) \) is invertible for all $k\in \mathbb{R}$, $n \in \mathbb{N}_0$, and $ \sup_{k, n} \left\| (\lambda - m_{\pm}(k, n))^{-1} \right\| < \infty $. The spectrum includes two branches of eigenvalues associated with the zero angular mode $n=0$:
\begin{equation*}
\lambda_{\pm}(k, 0) = i c k - \frac{1}{2} \left\{ D k^2 - C_\mathrm{m}^{-1} f'(\phi_1^\pm) + \eps \gamma 
+ \sqrt{ \left( D k^2 - C_\mathrm{m}^{-1} f'(\phi_1^\pm) - \eps \gamma \right)^2 - 4 \eps } \right\}\,.
\end{equation*}
Hence, $L_{\pm}$ are not sectorial and, consequently, the semigroups $e^{tL_{\pm}}$ are not analytic.

\begin{proof}[Proof of Proposition~\ref{prop:lin_semigroup}]
Using Lemma~\ref{lem:L_dissip} together with a classical perturbation result \cite[Theorem 3.1.1]{Pazy}, we conclude that the operator $L = L_{\pm} + V_{\pm}$ generates a strongly continuous semigroup $e^{tL}$ on $L^2$. 

To verify the estimate \eqref{eq:semigroupL} note that $-I+L$ is dissipative. This observation comes from the fact that
\begin{equation}
\label{eq:L-upperbound}
\Re \langle L v, v \rangle = \Re \langle (L_{\pm} + V_{\pm}) v, v \rangle \le \|v\|^2
\end{equation}
where we have used Lemma~\ref{lem:L_dissip} and the bound $ \langle V_{\pm} v, v \rangle \le \|v\|^2 $. Therefore, by Lemma~\ref{lem:dissipative} the semigroup $e^{tL}$ satisfies Eq.~\eqref{eq:semigroupL}, and the proof is complete.

\end{proof}

\subsection{Spectrum of the linearization $L$}
\label{subsec:spectrum_of_L}

In this subsection, we study the spectrum of $L$, which is crucial for constructing the projection operator $Q$ used in the decay analysis of Proposition~\ref{prop:LQ-decay} below.

\begin{lemma}\label{lem:ess}
	Let \( L \) be the operator defined in Eq.~\eqref{eq:L}. Its essential spectrum satisfies
	\[
	\sigma_{\mathrm{ess}}(L) \subset \left\{ \lambda \in \mathbb{C} \,\middle|\, \Re \lambda \le -\eta \right\},
	\]
	where \( \eta := \min\left\{ C_\mathrm{m}^{-1}|f'(\phi_1^\pm)|,\, \eps \gamma \right\} \).
\end{lemma}

\begin{proof}
	Since $L$ is a relatively compact perturbation of $L_{\pm}$ (see Subsection~\ref{subsec:semigroup}) we may use Weyl's essential spectrum theorem (Ch. 2, \cite{KPbook}). It follows that the essential spectrum remains unchanged under such perturbation, that is, $\sigma_{\mathrm{ess}}(L) = \sigma_{\mathrm{ess}}(L_{\pm})$. In Lemma~\ref{lem:L_dissip} we proved that $\eta + L_{\pm}$ is dissipative. Therefore, by Lemma ~\ref{lem:dissipative} the spectrum lies in the left half-plane, and the claim follows.
	
\end{proof}

\medskip

Next, we study the discrete spectrum of \( L \). We decompose functions on  \( \mathcal{S}_R \) via Fourier expansion in the angular variable \( \theta \). That is, for any \( v \in L^2(\mathcal{S}_R) \), we write
\[
v(z, \theta) = \sum_{n \in \mathbb{Z}} v_n(z)\, e^{i n \theta}.
\]
This decomposition allows us to reduce the operator \( L \) to a direct sum over decoupled mode operators $L = \bigoplus_{n \geq 0} L_n$, where each block \( L_n \) acts on the \( n \)-th Fourier mode (with \( n \geq 0 \) due to symmetry between \( \pm n \)) and is given by
\begin{equation}
\label{eq:Ln}
L_n :=
\begin{pmatrix}
D \left( \partial_z^2 - R^{-2}n^2 \right) + c \partial_z + C_\mathrm{m}^{-1}f'(\phi_1(z)) & -C_\mathrm{m}^{-1} \\
\eps & c \partial_z - \eps \gamma
\end{pmatrix}.
\end{equation}
The operators \( L_n \) describe the dynamics restricted to the invariant subspace corresponding to angular frequency \( n \).

In the following lemma we provide resolvent estimates for the positive modes.

\begin{lemma}\label{lem:Ln}
	Let \( L_n \) be defined as in Eq.~\eqref{eq:Ln}. If \( R \leq 1 \), then for all \( n > 0 \) we have:
	\begin{equation}\label{eq:spec_discr}
	\sigma_{\mathrm{d}}(L_n) \subset \left\{ \lambda \in \mathbb{C} \,\middle|\, \Re \lambda \le -\eta \right\},
	\end{equation}
	for $\eta := \min\left\{ C_\mathrm{m}^{-1}|f'(\phi_1^\pm)|,\, \eps \gamma \right\}$. The resolvent satisfies the estimate:
	\begin{equation}\label{eq:reso_Ln}
	\left\| (\lambda - L_n)^{-1} \right\| \le \frac{1}{\Re \lambda + \eta}, \qquad \text{for all } \Re \lambda > -\eta.
	\end{equation}
\end{lemma}

\begin{proof}
    It suffices to show that $\eta + L_n$ is dissipative for $n>0$. Then by Lemma~\ref{lem:dissipative} both Eq.~\eqref{eq:spec_discr} and the resolvent estimate of Eq.~\eqref{eq:reso_Ln} are satisfied. Indeed, using that $C_\mathrm{m}^{-1}(f'(\phi_1(z)) - f'(\phi_1^{\pm})) \le 1$ we have
	\begin{align*}
	\Re\, \langle L_n v, v \rangle
	&= \int_{\mathbb{R}} \left( D\, (\partial_z^2 v_1) \bar{v}_1 + \left( C_\mathrm{m}^{-1} f'(\phi_1(z)) - D\, n^2 R^{-2} \right) |v_1|^2 - \gamma |v_2|^2 \right) R\, dz \\
	&\le -\eta \|v\|^2,
	\end{align*}
	provided that $n>0$ and $R\leq 1$.
	
\end{proof}

We now focus on the discrete spectrum of the mode $n = 0$, which carries the neutral eigenvalue.

\begin{lemma}\label{lem:L0}
	Let \( L_0 \) denote the operator \( L_n \) defined in Eq.~\eqref{eq:Ln} with \( n = 0 \). Then, for sufficiently small \( \eps > 0 \), there exists \( \beta \gtrsim \eps \) such that
	\[
	\sigma_{\mathrm{d}}(L_0) \subset \left\{ 0 \right\} \cup \left\{ \lambda \in \mathbb{C} \,\middle|\, \Re \lambda \le -\beta \right\},
	\]
	and the eigenvalue at zero is simple for both $L_0$ and its adjoint $L_0^*$.
\end{lemma}

\begin{proof}
	Let $\tau := -\partial_z \Phi$ denote the tangent vector to $\mathcal{M}$ at $\Phi$. As shown earlier, the operator $L$ satisfies $L \tau = 0$ due to translation invariance, so $0$ is an eigenvalue of $L$ with eigenfunction $\tau$. It remains to show that $0$ is a simple eigenvalue, and that there are no other eigenvalues with non-negative real part. 
    Under the assumption that $\eps$ is positive and sufficiently small, the argument follows in a manner similar to that of the proof of Lemma 3.10 in \cite{Talidou_21} and here is omitted.
	
\end{proof}

Therefore, by Lemmas~\ref{lem:ess},~\ref{lem:Ln}, and~\ref{lem:L0} we conclude that
\begin{equation}\label{eq:spectrum}
\sigma(L) \subset \{0\} \cup \left\{ \lambda \in \mathbb{C} \,\middle|\, \Re \lambda \le -\eta' \right\},
\end{equation}
where $\eta' := \min\left\{ C_\mathrm{m}^{-1} |f'(\phi_1^\pm)|,\ \beta,\ \eps \gamma \right\}$, and $\beta$ is a spectral gap determined in Lemma~\ref{lem:L0}.

\subsection{The spectral projection}
\label{sec:projection}

Since the spectrum of $L$ contains the zero eigenvalue, which is isolated, we define the associated Riesz projection using the standard contour integral representation:
\[
P = \frac{1}{2\pi i} \oint_{\Gamma_0} (\lambda - L)^{-1}\, d\lambda,
\]
where $\Gamma_0 \subset \mathbb{C}$ is a positively oriented closed contour enclosing the origin and no other spectral values of $ L $ \cite[Chapter 2]{KPbook}. The projection $P$ is a bounded linear operator that commutes with $L$ and its semigroup. The complementary projection is defined by $Q := I-P$, and is the projection onto $\operatorname{Ran}(L)$. $Q$ also commutes with $L$ and $ e^{tL} $.
The decomposition via the Riesz projection \(P + Q = I\) is crucial for isolating the neutral mode due to translational symmetry. 
%You could maybe emphasise that semigroup decay only applies on \(\mathrm{Ran}(Q)\), making this spectral splitting essential

\begin{lemma}\label{lem:P}
	Let $\tau = -\p_z \Phi$ be the tangent vector on $\cM$ at a fixed point $\Phi$, and $\tau^*$ the eigenfunction of the adjoint $L^*$ corresponding to the zero eigenvalue, normalized such that $\langle \tau, \tau^* \rangle = 1$. Then the Riesz projection is given explicitly by
	\[
	P v = \langle v, \tau^* \rangle\, \tau, \qquad v \in L^2.
	\]
\end{lemma}

\begin{proof}
	Since both $L$ and $L^*$ possess a simple eigenvalue at zero (as shown in Lemma~\ref{lem:L0}), the eigenfunction $\tau^*$ is uniquely determined up to scalar normalization. Moreover, because $\tau \notin \operatorname{Ran}(L)$, is not orthogonal to $\tau^*$. The normalization ensures that $P \tau = \tau$ and $P L v = 0$ for all $v \in H^{2,1}$. To verify the form, observe that $\langle \tau, \tau^* \rangle\, \tau = \tau$ and $\langle L v, \tau^* \rangle\, \tau = \langle v, L^* \tau^* \rangle\, \tau = 0$, since \( L^* \tau^* = 0 \). This spectral projection isolates the neutral mode associated with translational symmetry, a typical feature in traveling wave problems.
	
\end{proof}

\subsection{Resolvent estimate for the zero Fourier mode}
\label{subsec:resolvent-estimates}

Next, we establish the boundedness of the resolvent operator $(\lambda - L_0)^{-1}$. Let $P$ and $Q$ be the spectral projections onto $\operatorname{Ker}(L)$ and $\operatorname{Ran}(L)$, respectively, and denote by $Q_0$ the restriction of $Q$ to the zero Fourier mode subspace.

\begin{lemma}\label{lem:Res0}
	Consider the operator $L_0$ defined as in Eq.~\eqref{eq:Ln} with $n=0$. For every $\mu < \eta':= \min\{C_\mathrm{m}^{-1} |f'(\phi_1^\pm)| , \beta, \eps \gamma\}$ there exists a constant $C>0$ such that
	\begin{equation}
	\| (\lambda - L_0)^{-1} Q_0 \| \le C, \qquad \forall\, \lambda \in \mathbb{C},\ \Re\,\lambda \ge -\mu.
	\end{equation}
\end{lemma}

\begin{proof}

We divide the region $\Re\,\lambda \ge -\mu$ into three distinct subdomains:
\begin{align*}
S_1 &:= \left\{ \lambda \in \mathbb{C} \mid -\mu \le \Re\,\lambda \le 2,\ |\Im\,\lambda| \le N \right\}, \\
S_2 &:= \left\{ \lambda \in \mathbb{C} \mid \Re\,\lambda \ge 2 \right\}, \\
S_3 &:= \left\{ \lambda \in \mathbb{C} \mid -\mu \le \Re\,\lambda \le 2,\ |\Im\,\lambda| \ge N \right\},
\end{align*}
for some threshold \( N > 0 \) to be chosen later, and estimate the resolvent operator of $L_0$ in each of these.

\paragraph{Estimate on \( S_1 \):}
Since the spectrum of \( L_0 \) intersects \( S_1 \) only at the simple eigenvalue at zero (see Lemma~\ref{lem:L0}), the operator \( (\lambda - L_0)^{-1} \) is analytic on $\operatorname{Ran}(Q_0)$ throughout \( S_1 \). The compactness of \( S_1 \) and analyticity of the resolvent then imply that for any $N>0$
\begin{equation}\label{eq:estimate_S1}
\sup_{\lambda \in S_1} \| (\lambda - L_0)^{-1} Q_0 \| < \infty.
\end{equation}

\paragraph{Estimate on \( S_2 \):} Due to dissipativity estimate from \eqref{eq:L-upperbound}, we have $\Re\,\langle L_0 v, v \rangle \le \|v\|^2$. Applying Lemma~\ref{lem:dissipative}, it follows that for \( \Re\lambda > 1 \), the resolvent satisfies
\begin{equation}\label{eq:estimate_S2}
\|(\lambda - L_0)^{-1} \| \le \frac{1}{\Re\lambda - 1} \le 1\,,
\end{equation}
for all $\lambda \in S_2$. The second inequality follows because $\Re\lambda \ge 2$ in $S_2$.

\paragraph{Estimate on \( S_3 \):} To provide resolvent estimates on this subdomain we lean on the proof of Lemma 3.15 in \cite{Talidou_21}. Specifically, we prove the following explicit estimate.

\begin{lemma}
	\label{lem:S3}
	There exists $N > 0$ sufficiently large such that
	\begin{equation}\label{eq:estimate_S3}
	\sup_{\lambda \in S_3} \| (\lambda - L_0)^{-1} \| \le \frac{2}{\eta' - \mu}
	\end{equation}
	where $\eta' := \min\left\{ C_\mathrm{m}^{-1} |f'(\phi_1^\pm)|,\ \beta,\ \eps \gamma \right\} > \mu$.
\end{lemma}

\begin{proof}
	We view $L_0$ as a perturbation of the constant-coefficient operator $L_{0}^{\pm}$:
	\begin{equation*}
	L_{0}^{\pm} := 
	\begin{pmatrix}
	D \p_z^2 + c \p_z + C_\mathrm{m}^{-1} f'(\phi_1^\pm) & -1 \\
	\eps & c \p_z - \eps \gamma
	\end{pmatrix},
	\end{equation*}
	and define \( V := L_0 - L_{0}^{\pm} \). The key estimate follows from showing that for sufficiently large \( N \), the operator norm \( \| (\lambda - L_{0}^{\pm})^{-1} V \| < 1/2 \) uniformly on \( S_3 \). This is established by direct analysis of $(\lambda - m_{\pm}(k,0))^{-1}$, where $m_{\pm}(k,0)$ is the Fourier multiplier associated with $L_{0}^{\pm}$, as shown in Eq.~\eqref{eq:multOp}. Through detailed asymptotic estimates on each matrix entry using Cramer's rule, and separating real and imaginary parts, one obtains decay of the operator norm as \( |\Im\,\lambda| \to \infty \), ensuring smallness of the perturbation term. 
	
	The desired estimate follows by the resolvent identity:
	\begin{equation*}
	(\lambda - L_0)^{-1} = \left( I - (\lambda - L_{0}^{\pm})^{-1} V \right)^{-1} (\lambda - L_{0}^{\pm})^{-1},
	\end{equation*}
	and the dissipativity of $\eta' + L_{0}^{\pm}$.
		
	For a more detailed proof we refer the reader to \cite[Lemma 3.15]{Talidou_21}.

\end{proof}

\medskip

\noindent The uniform boundedness of the resolvent $(\lambda - L_0)^{-1} Q_0$ for $\Re\lambda \ge -\mu$ now follows directly from the estimates \eqref{eq:estimate_S1}-\eqref{eq:estimate_S3}, completing the proof of Lemma~\ref{lem:Res0}.

\end{proof}

\subsection{Exponential decay of the linear semigroup}
\label{subsec:LQ-decay}

We are now able to prove the exponential decay of the semigroup $e^{tL}$ on the subspace orthogonal to the kernel of $L$.

\begin{proposition}
	\label{prop:LQ-decay} 
	Let $L$ be the linearized operator defined in Eq.~\eqref{eq:L}. Assume that \( \eps > 0 \) is sufficiently small and that the cylinder radius satisfies \( 0 < R \leq 1 \). Then there exists a constant \( \eta' > 0 \) such that the semigroup \( e^{tL} \), restricted to $\operatorname{Ran}(Q)$, satisfies the exponential decay estimate
	\begin{equation} \label{eq:LQ-decay}
	\big\| e^{tL} Q \big\|_{2,1} \leq C e^{-\eta' t}, \qquad \text{for all } t \geq 0,
	\end{equation}
	for some constant \( C > 0 \).
\end{proposition}

\begin{proof}
    Let $P$ and $Q$ be the spectral projections as defined in Subsection~\ref{sec:projection}, and choose $\mu >0$ sufficiently small so that $\mu < \eta'$ where recall that $ \eta' := \min\left\{ C_\mathrm{m}^{-1} |f'(\phi_1^\pm)|,\, \beta,\, \eps \gamma \right\}$, and $\beta$ is the spectral gap in Lemma~\ref{lem:L0}. 
		
    The operator $L$ decomposes into Fourier modes $L = \bigoplus_{n \ge 0} L_n$, and hence
    \begin{equation}\label{eq:reso_L_bound}
	\| (\lambda - L)^{-1} Q \| = \sup_{n \ge 0} \| (\lambda - L_n)^{-1} Q_n \| < C,
    \end{equation}
    for some constant $C>0$. The uniform bound of the resolvent follows from Lemma~\ref{lem:Ln} (for $n>0$), and Lemma~\ref{lem:Res0} (for $n=0$). Since $L$ commutes with $Q$, by the equivalence relation of Lemma~\ref{lem:L-domain}, the resolvent estimate \eqref{eq:reso_L_bound} also holds in the $\| \cdot \|_{2,1}$. 
	
    For $n > 0$, Lemma~\ref{lem:Ln} provides uniform bounds on the resolvent, while Lemma \ref{lem:Res0} ensures the same for $n = 0$. Therefore, the resolvent of $L$ restricted to $\operatorname{Ran}(Q)$ is bounded for $\Re\,\lambda \ge -\mu$, and we may apply \cite[Corollary 4 (Pr\"uss' theorem)]{Pruss}  to deduce that
    \begin{equation*}
        \| e^{tL} Q \| \le C e^{-\mu t}, \qquad \forall\, t \ge 0\,,
    \end{equation*}
    for $C>0$. Since $L$ commutes with $e^{tL}$,
    \begin{equation*}
	\| e^{tL} Q \|_{2,1} \le \| L e^{tL} Q u \| + \| e^{tL} Q u \| \le C e^{-\mu t} ( \| L u \| +      \| u \| ) \le C e^{-\mu t} \| u \|_{2,1}\,.
    \end{equation*}
    The first and last inequalities follow from the equivalence of the graph norm of $L$ with the $H^{2,1}$-norm (cf. Lemma~\ref{lem:L-domain}), completing the proof of Proposition~\ref{prop:LQ-decay}.

\end{proof}

\noindent %As a direct implication of this result, the full semigroup \( e^{tL} \) remains uniformly bounded in the \( H^{2,1} \)-norm for all time. This follows from the triangle inequality:
%\[
%\|e^{tL}\|_{2,1} \leq \|e^{tL}P\|_{2,1} + \|e^{tL}Q\|_{2,1} \leq \|P\|_{2,1} + C \|Q\|_{2,1}, \qquad \text{for all } t \geq 0,
%\]
%since \( e^{tL} \) acts as the identity on the range of \( P \) and decays exponentially on the complementary subspace projected by \( Q \).

\subsection{Proof of Theorem~\ref{thm:nonli-stab}}
As mentioned above, the proof of Theorem~\ref{thm:nonli-stab} relies on the application of Theorem 4.3.5 from \cite{KPbook} (also restated in Appendix~\ref{app:stability_thm}). It therefore remains to verify that all assumptions of the theorem are satisfied. Assumption (a) follows from Eq.~\eqref{eq:spectrum}, and (b) from Eq.~\eqref{eq:span}. The third assumption is established in Proposition~\ref{prop:LQ-decay}, specifically in Eq.~\eqref{eq:reso_L_bound}. We also note that Proposition~\ref{prop:LQ-decay} goes a step further by providing, in addition to resolvent estimates, bounds on the semigroup $e^{tL}$, which will be used later in the proof of Theorem~\ref{thm:persist}. Moreover, there exist positive constants $M_1$, $M_2$ and $m$ such that
\begin{equation}\label{eq:nablaF}
	\begin{split}
	\| \nabla_uF(u) - \nabla_uF(v) \|_{2,1} &\leq \| C_\mathrm{m}^{-1} \left(f'(u_1) - f'(v_1) \right) \|_{H^2} \\
	&\leq \| C_\mathrm{m}^{-1} \left( 3(u_1^2 + v_1^2) + 2(\alpha + 1)(u_1 - v_1) \right) \|_{H^2} \\
	&\leq M_1 \| u-v \|_{2,1}
	\end{split}
\end{equation}
and
\begin{equation}\label{eq:nonlin}
\begin{split}
\| N(v) \|_{2,1} &= \| C_\mathrm{m}^{-1} v_1^2 (-v_1 -3\phi_1 + \alpha + 1) \|_{H^2} \\
&\leq M_2 \| v \|_{2,1}^2
\end{split}
\end{equation}
for all $\| u_1 \|_{H^2} \leq m$ and $\| v_1 \|_{H^2} \leq m$. The estimates \eqref{eq:nablaF} and \eqref{eq:nonlin} verify assumptions (d) and (e), respectively. Therefore, all conditions of \cite[Theorem 4.3.5]{KPbook} are satisfied, which completes the proof of Theorem~\ref{thm:nonli-stab}.

\section{Front-like solutions on warped cylinders}
\label{sec:warped}
In this section, we consider warped cylinders $\cS = \cS_{\rho}$ as defined in Eq.~\eqref{eq:def-Srho}, and prove Theorem~\ref{thm:persist}. Recall that in the proof of Theorem~\ref{thm:nonli-stab}, fluctuations around the traveling front were controlled by comparing the full system to the one linearized around the front in a moving frame. In the present setting, however, the spatial variation of the radius induces a time-dependent perturbation in the principal part of the linearized operator, as the operator $\Del_{\cS_{\rho(x)}}$ in Eq.~\eqref{eq:Delta-rho} becomes $\Del_{\cS_{\rho(z+ct)}}$. This results in a time-dependent evolution system that is challenging to analyze. To address this, we adopt the idea from \cite{Talidou_21} and treat the flow on $\cS_{\rho}$ in the static frame as a perturbation of the flow on the standard cylinder $\cS_{R}$. We first estimate the difference between the flows linearized at $u=0$, and then, by applying Gronwall's inequality, extend these estimates to the full nonlinear flow.

Due to the geometry of the surface of $\cS_{\rho}$, the inner product on $L^2_{\rho} \equiv L^2(\cS_\rho)$ is defined by the surface integral
\begin{equation}\label{eq:inner-rho}
\langle u,w \rangle_\rho :=\int_{\cS_{\rho}} 
(u_1\bar{w}_1+\eps^{-1}u_2\bar{w}_2)\, d\mu_\rho\,,
\end{equation}
with the corresponding norm $\| \cdot \|_\rho$. In Eq.~\eqref{eq:inner-rho}, $d\mu_\rho = \sqrt{g}\, d\theta dx$ is the Riemannian area element with density $g= \rho^2(1+\rho'^2)$.
Moreover, define the mixed Sobolev spaces
\begin{align}\label{eq:Hkl}
H^{2k,\ell}_\rho
:=\left\{u\in L^2
\ \big\vert\  
(\Delta_{\cS_\rho})^k u_1 \in L^2_\rho,
(\p_x)^\ell u_2\in L^2_\rho \right\}
\end{align}
for $k,\ell = \{0,1\}$,
with norms 
\begin{equation}\label{eq:norms_k_l}
\|u\|_{2k,\ell;\rho}:=\sum_{0\le i\le k}
\|(\Delta_{\cS_\rho})^i u_1\|_\rho + 
\eps^{-1}\sum_{0\le j\le \ell}
\|\p_x^j u_2\|_\rho\,.
\end{equation}
For $k=1, \ell=0$, the space $H^{2, 0}$ agrees with 
the corresponding Sobolev space $H^{2}\times L^2$,
and $H^{0,0}=L^2$. However,
for $\ell=1$, since $H^{2k,1}$
places no condition on $\p_\theta u_2$,
the space $H^{0,1}$ properly contains
$L^2\times H^1$, and $H^{2,1}$ properly contains
$H^{2}\times H^1$. 
On the standard cylinder~$\cS_R$,
Eq.~\eqref{eq:Hkl} with $k=\ell=1$ coincides
with the definition of $H^{2,1}$ in Eq.~\eqref{eq:H21}.

\subsection{Linear dynamics}

Let $F_{\rho}(u)$ denote the right hand side of Eq.~\eqref{eq:FHN} with $\cS = \cS_{\rho}$. The traveling front $\Phi$ is not a solution to that equation, i.e., $F_{\rho}(\Phi) \neq 0$. Since there is no advantage if we Taylor expand Eq.~\eqref{eq:FHN} around $\Phi$ we choose to expand around the zero solution. The linearized operator, denoted by $A_{\rho}$, is defined by the G\^ateaux derivative:
\begin{equation} \label{eq:A} 
A_\rho:= dF_\rho(0)
= \begin{pmatrix} C_\mathrm{m}^{-1} (\Delta_{\cS_{\rho}} + f'(0)) & -1 \\
\eps & - \eps\gamma 
\end{pmatrix}\,,
\end{equation}
where $\Delta_{\cS_{\rho}}$ is defined in Eq.~\eqref{eq:Delta-rho}.

\subsubsection{Spectral properties of $A_{\rho}$}
\label{sec:properties_A_rho}

The linearization $A_{\rho}$ defines a closed linear operator on $L^2_{\rho}$. The domain of $A_{\rho}$ is $H^{2,0}_{\rho}$, and its graph norm is equivalent to $\| \cdot \|_{2,0;\rho}$. The graph norm of $A_{\rho}$ is also equivalent to $\| A_{\rho}u \|_{\rho}$:
\begin{equation}
\label{eq:domain-A}
\|u\|_{2,0;\rho}\lesssim \|A_\rho u\|_\rho
\lesssim \|u\|_{2,0;\rho}\,.
\end{equation}
Equation \eqref{eq:domain-A} follows from the fact that $A_{\rho}$ is dissipative. Specifically, 
\begin{align*}
\Re\,\langle A_\rho u, u\rangle_\rho &= \langle C_\mathrm{m}^{-1} (\Delta_{\cS_{\rho}} + f'(0)) u_1, u_1 \rangle - \gamma \| u_2 \|^2 \\ 
&\le -\nu \|u\|_\rho^2\,,
\quad u\in H^{2,1}_{\rho}\,,
\end{align*}
where $\nu := \min\left\{ C_\mathrm{m}^{-1} \alpha, \eps \gamma \right\}$ and $\langle \cdot, \cdot \rangle_{\rho}$ is the inner product defined in Eq.~\eqref{eq:inner-rho}. By that same property and Lemma~\ref{lem:dissipative}, $A_{\rho}$ generates a strongly continuous semigroup of contractions $e^{t A_{\rho}}$, that decays exponentially on $L^2_{\rho}$ and has $H_{\rho}^{2,0}$ as an invariant subspace.

Since for the comparison to the solutions on $\cS_{R}$ we will work in the subspace $H_R^{2,1}$ (that was used for Theorem~\ref{thm:nonli-stab}), we restrict $A_{\rho}$ to the intermediate subspace $H^{0,1}_{\rho}$ and provide spectral estimates. Specifically, the operator $A_{\rho}$ maps $H^{2,1}_{\rho}$ bijectively onto $H_{\rho}^{0,1}$, and
\begin{equation}\label{eq:domain-A-01}
\| u \|_{2,1;\rho} 
\lesssim \| A_\rho u \|_{0,1;\rho}
\lesssim \| u \|_{2,1;\rho}
\end{equation}
for all $u \in H^{2,1}_{\rho}$ (see Lemma 5.2 in \cite{Talidou_21}). As a consequence of Eq.~\eqref{eq:domain-A-01}, for every bounded linear operator $B$ on $H^{0,1}_{\rho}$ that commutes with $A_{\rho}$ it holds that
\begin{equation*}
\|B\|_{0,1;\rho}\lesssim\|B\|_{2,1;\rho}\lesssim \|B\|_{0,1;\rho}\,.
\end{equation*}
Next we study the spectrum of $A_{\rho}$ on $H^{0,1}_{\rho}$ and provide estimates on its semigroup and resolvent.

\begin{lemma}\label{lem:sectorial} Let $\rho\in C^2$ be a real-valued function that is bounded and bounded away from zero. Then, $A_\rho$ generates an analytic semigroup $e^{tA_\rho}$ on $H^{2,1}$. The spectrum of $A_\rho$ on $H^{0,1}_{\rho}$ is contained in the truncated sector
\begin{equation}\label{eq:sector}
\Sigma:= \left\{ \lambda\in \bC \ \big\vert \ 
\Re\,\lambda \le 
-\nu \min\left\{1-\frac{C_\mathrm{m} r_{int}}{4}, \eps^{-\frac12}\, | \Im\, \lambda|\right\} \right\} \,,
\end{equation}
where $\nu:=\min\left\{C_\mathrm{m}^{-1} \alpha, \eps\gamma\right\}$.
\end{lemma}

\begin{proof}
	For $s>0$, define the inner product 
	\begin{equation*}
	\langle u, w \rangle_s :=
	\langle u_1, w_1 \rangle_{\rho} 
	+ \eps^{-1}\langle u_2, w_2 \rangle_{\rho} + s\eps^{-1}
	\langle \p_x u_2, \p_x w_2\rangle_{\rho}\,.
	\end{equation*}
	From Eq.~\eqref{eq:Hkl}, the corresponding norm $\|\cdot\|_s$ is equivalent to $\|\cdot\|_{0,1}$:
	\begin{equation}\label{eq:equiv_s}
	\min\{1,s^\frac12\} \|u\|_{0,1;\rho}\le \|u\|_s
	\le \max\{1,s^{\frac12}\}\|u\|_{0,1;\rho}\,,\qquad (s>0)\,. 
	\end{equation}
	For the real part of $\langle A_{\rho}u, u \rangle_s$ we have: 
	\begin{align*}
	\notag \Re\, 
	\langle A_\rho u, u \rangle_s
	&= \langle C_\mathrm{m}^{-1} \Delta_{\cS_\rho} u_1, u_1 \rangle_{\rho} - C_\mathrm{m}^{-1}\alpha \|u_1\|^2_{\rho} -\gamma \|u_2\|^2_{\rho} + s\Re\,  
	\langle \p_x u_1-\gamma\p_xu_2, \p_x u_2\rangle_{\rho}\\
	\notag
	&\hskip -1cm\le -\frac{\pi}{C_m r_{int}} \frac{\sup|\rho|^2}{1+\sup|\rho'|^2}\|\p_x u_1\|^2_{\rho} - \nu \|u\|^2_{\rho} + s \|\p_x u_1\|_{\rho} \|\p_x u_2\|_{\rho}
	-\gamma\|\p_x u_2\|^2_{\rho}\,.  
	\end{align*}
	In the last line of the above computation we have integrated the Laplacian by parts. By completing the square we have
	\begin{equation*}
	\Re\, \langle A_\rho u, u \rangle_s \leq -\nu\|u\|_s^2 + \frac{s^2}4 \frac{C_m r_{int}}{\pi} \frac{1 + \sup|\rho'|^2}{\sup|\rho|^2}\|\partial_x u_2\|^2\,.
	\end{equation*}
	Choosing $s \leq \frac{\sup|\rho|^2 \pi \nu}{\eps} $ the above yields
	\begin{equation}\label{eq:Re_A_rho}
	\Re\, \langle A_\rho u,u\rangle_s\le -\nu\left(1-\frac{C_m r_{int}}{4}\right)\|u\|_s^2\,.   
	\end{equation}
	Since $A_{\rho}$ is dissipative, Lemma~\ref{lem:dissipative} implies that $A_{\rho}$ generates a strongly continuous semigroup of contractions with respect to $\| \cdot \|_s$. The spectrum of $A_{\rho}$ is contained in each of the half-planes $\{\Re\, \lambda\le-\nu\left(1-\frac{C_m r_{int}}{4}\right)\}$, and hence in their intersection.
	
	For the imaginary part of $\langle A_{\rho}u, u \rangle_s$ we have:
	\begin{align*}
	\Im\, \langle A_\rho u,u\rangle_s 
	& = (1 + C_\mathrm{m}^{-1}) \Im\, \langle u_1, u_2\rangle_{\rho}
	+ s\Im\, \langle \p_x u_1,\p_x u_2\rangle_{\rho}\\
	&\le \sqrt{\eps}\|u\|^2_{\rho} + s\|\p_x u_1\|_{\rho}
	\|\p_xu_2\|_{\rho}\,.
	\end{align*}
	Comparing this with Eq.~\eqref{eq:Re_A_rho}, we observe that for $s>0$ sufficiently small
	\begin{equation*}
	\Re\,\langle A_\rho u, u\rangle_s \le -\nu\eps^{-\frac12}
	|\Im\, \langle A_\rho u, u\rangle_s|\,.  
	\end{equation*}
	Since the resolvent set of $A_\rho$ contains 0, by \cite[Theorem~1.3.9]{Pazy}, it contains the entire complement of the sector $\{\Re\,\lambda\le 
	- \nu \eps^{-\frac12}|\Im\, \lambda|\}$. In summary, for $s>0$ sufficiently small, the
	numerical range of $A_\rho$ with respect to $\langle\cdot,\cdot\rangle_s$ lies in \eqref{eq:sector}. 
	
\end{proof}

\noindent By Lemma~\ref{lem:sectorial} and (\cite{Pazy}, Theorem 1.3.9) the resolvent operator satisfies
\begin{equation*}
\|(\lambda-A_{\rho})^{-1}\|_s \le \left(\inf_{z\in \Sigma}
\|\lambda-z\|_s\right)^{-1}\,,\qquad (\lambda\not\in \Sigma)\,.
\end{equation*}
For $s= \gamma(1+\sup|\rho'|^2)^{-1}$ and the equivalence relation \eqref{eq:equiv_s} we have
\begin{equation}\label{eq:sectorial}
\|(\lambda-A_{\rho})^{-1}\|_{0,1;\rho}\le 
C(1+\sup|\rho'|) \min\left\{1,|\lambda|^{-1}\right\}\,,
\end{equation}
for all $\lambda$ with $\Re \lambda \geq -\frac{4-C_m\, r_{int}}{8} \nu \min\left\{ 1, \eps^{-\tfrac{1}{2}} |\Im \lambda| \right\}$ and some constant $C$. Since the resolvent commutes with $A_{\rho}$, by Eq.~\eqref{eq:domain-A-01} the estimate \eqref{eq:sectorial} holds also for the norm $\| \cdot \|_{2,1;\rho}$ (with a suitably adjusted constant). By Lemma~\ref{lem:dissipative} and Eq.~\eqref{eq:sectorial} we have the following estimate for the semigroup generated by $A_{\rho}$:
\begin{equation}
\sup_{t>0} \| e^{tA_{\rho}} \|_{2,1;\rho} \leq C(1+\sup|\rho'|) e^{-\nu t}\,.
\end{equation}

%%%%%%%%% 
\subsubsection{Perturbation estimate for the semigroup}
\label{subsec:compare}

The main result of this subsection is the following lemma:

\begin{lemma}\label{lemma:difference}
    Let $\alpha$, $\gamma$, $\eps$, $C_\mathrm{m}$ and $r_\mathrm{int}$ be positive constants. There exists a constant $C$ such that, if  $0<R\le 1$ and $\delta:=R^{-1}\|\rho-R\|_{C^2}\le c$, for some $c>0$, then the semigroup generated by $A_\rho$ on $H^{2,1}$ satisfies the estimate
    \begin{equation}\label{eq:semigroup_dif}
        \| e^{t A_R} - e^{t A_{\rho}} \|_{2,1} \leq C\delta
	(1 + \log (t^{-1} + 1))
    \end{equation}
    for all $t> 0$. 
\end{lemma}

To prove this lemma, we compare solutions of Eq.~\eqref{eq:FHN} on $\cS_{\rho}$ with solutions on $\cS_R$. In Eq.~\eqref{eq:diffeo} we have defined the diffeomorphism $\psi_\rho:\cS_R\to\cS_\rho$ which allows us to identify functions $u$ on $\cS_{\rho}$ with functions on $\cS_R$. We next show that $u\mapsto u\circ\psi_\rho$ induces a bounded linear transformation from $H^{2,1}(\cS_\rho)$ to $H^{2,1}(\cS_{R})$.

\begin{lemma}\label{lem:21-equivalent} 
    Let $\alpha$, $\gamma$, $\eps$, $C_\mathrm{m}$ and $r_\mathrm{int}$ be fixed positive 
    constants, $0 < R \leq 1$ and $\rho \in C^2$. There exists a constant $c>0$ such that if           $\delta:=R^{-1}\|\rho-R\|_{C^2}\le c$, then
    \begin{equation*}
	\|u\|_{2k,\ell} \lesssim \|u\|_{2k,\ell;\rho} \lesssim \|u\|_{2k,\ell}\,
	\qquad (k,\ell\in\{0,1\})\,.
    \end{equation*}
\end{lemma}

\begin{proof}
	From the definition of the norms $\| \cdot \|_{2k, \ell;\rho}$, $k, \ell\in \{0,1\}$, in Eq.~\eqref{eq:norms_k_l} we need lower and upper bounds for $\| w \|_{\rho}$, $\| \p_x w \|_{\rho}$ and $\| \Delta_{\cS_{\rho}} w \|_{\rho}$.
	
	By Eq.~\eqref{eq:inner-rho} and the pointwise bound $| \sqrt{g(x)} - R | \leq 2\delta$ we have that
	\begin{equation}\label{eq:norm_w_rho}
	\big| \|w\|_{\rho}^2 - \|w\|^2\big| \leq c\delta \|w\|^2\,.
	\end{equation}
	Solving for the norm on $L^2_{\rho}$ we obtain
	\begin{equation}\label{eq:norm_w_rho_2}
	\| w \| \lesssim \| w \|_{\rho} \lesssim \| w \|\,.
	\end{equation}
	In the same manner we obtain bounds for $\| \p_x w \|_{\rho}$.
	
	Next, recall that $A_{\rho}$ is given by Eq.~\eqref{eq:A}, and the corresponding operator $A_R$ is defined as
	\begin{equation} \label{eq:A_R} 
	A_R	= 
	\begin{pmatrix} 
	D\Delta_{\cS_R} + C_\mathrm{m}^{-1} f'(0) & -1 \\
	\eps & - \eps\gamma 
	\end{pmatrix}\,.
	\end{equation} 
	To bound the norm $\| \Delta_{\cS_{\rho}} w \|_{\rho}$ we consider the difference
	\begin{equation*}
	C_\mathrm{m}^{-1} \Del_{\cS_{\rho}} - D \Delta_{\cS_R} = -\nu_1(x) \p_x^2 + \nu_2(x) \p_x\,,
	\end{equation*}
	where
	\begin{align*}
	\nu_1(x) &:= -\frac{\pi R^2}{C_m r_{int}} + \frac{\pi}{C_m r_{int}} \frac{\rho'(3(1+(\rho')^2) - \rho \rho'')}{(1+(\rho')^2)^2} \\
	\nu_2(x) &:= \frac{\pi}{C_m r_{int}} \frac{\rho^2}{1+(\rho')^2}\,.
	\end{align*}
	These two coefficients are controlled by $\delta$. Since $\delta\le1$, $R\le 1$, and
	$\|\partial_x w\|\le \frac12(\|D\Delta_{\cS_R} w\| +\|w\|)$, by the triangle inequality we have
	\begin{equation}
	\label{eq:rho-R-4}
	\left\| (C_\mathrm{m}^{-1} \Delta_{\cS_\rho} - D \Delta_{\cS_R}) w\right\|
	\le  c\delta \left( \|D\Delta_{\cS_R}w\| + \|w\|\right)\,,
	\end{equation}
	where $c>0$ is independent on $\rho$ and $R$. Using the triangle inequality together with Eq.~\eqref{eq:norm_w_rho},
	we arrive at
	\begin{equation*}
	\bigl| \|C_\mathrm{m}^{-1}\Delta_{\cS_\rho}w\|_\rho - \|D\Delta_{\cS_R}w\| \bigr|
	\leq c\delta \left( \|D\Delta_{\cS_R}w\| + \|w\|\right)\,.
	\end{equation*}
	Solving with respect to the $L^2_{\rho}$-norm we obtain bounds for $\| \Del_{\cS_{\rho}} w \|_{\rho}$.
	
\end{proof}

For the rest of this section, we identify the spaces $H^{2k,\ell}(\cS_\rho)$ with $H^{2k,\ell}(\cS_R)$, and use the standard norms $\|\cdot\|_{2k,\ell}$. We now prove Lemma~\ref{lemma:difference}.

\paragraph{Proof of Lemma~\ref{lemma:difference}.} Since $A_{\rho}$ is sectorial (see Lemma~\ref{lem:sectorial}), the semigroup $e^{tA_{\rho}}$ is given by the contour integral
\begin{equation*}
    e^{tA_\rho} = \frac{1}{2\pi i}
    \oint_{\Gamma} e^{\lambda t} (\lambda-A_\rho)^{-1}\,d\lambda\,,
\end{equation*}
where $\Gamma$ is the contour consisting of the two half-lines $Re \lambda = -\frac{4-C_m r_{int}}{8} \nu \eps^{-\frac{1}{2}} |Im \lambda|$ traversed counterclockwise, enclosing the spectrum of $A_{\rho}$. A similar contour is defined for the operator $A_R$. Parametrizing $\Gamma$ by $\lambda(s) = -|s|+ i \frac{8}{4-C_m r_{int}} \nu^{-1}\sqrt\eps$, we observe that for each $t>0$, the above integral converges absolutely with respect to the operator norm on $H^{2,1}_R$.

Next, we estimate the difference of $e^{tA_{\rho}}$ from $e^{tA_R}$:
\begin{equation}\label{eq:difference_semigroups}
    \begin{split}
        \|e^{tA_\rho}-e^{tA_R}\|_{2,1}
        &= \left\|\frac{1}{2\pi i}
        \oint_{\Gamma} e^{\lambda t} 
        \bigl( (\lambda-A_\rho)^{-1}- (\lambda-A_R)^{-1}\bigr)\,d\lambda\right\|_{2,1}\\
        &\le \frac{1}{2\pi}
        \int_{\mathbb{R}} e^{t\Re \lambda(s)}
        \|(\lambda\!-\!A_\rho)^{-1} -(\lambda\!-\!A_R)^{-1}\|_{2,1}\, 
        |\lambda'(s)|\, ds\,.    
    \end{split}
\end{equation}
It remains to estimate the difference of the resolvents. By the second resolvent identity we have
\begin{equation}\label{eq:reso_id}
    (\lambda\!-\! A_{\rho})^{-1} - (\lambda\!-\! A_{R})^{-1} =
    (\lambda\!-\! A_R)^{-1} W(\lambda\!-\!A_\rho)^{-1}
\end{equation}
where
$$
W := A_\rho-A_R = 
\begin{pmatrix}
    C_\mathrm{m}^{-1}\Delta_{\cS_\rho} -D\Delta_{\cS_R} & 0 \\
    0 & 0
\end{pmatrix}\,.
$$
By Eq.~\eqref{eq:norm_w_rho_2} the operator $W:H^{2,1} \to H^{0,1}$ satisfies
$$
\|W u\|_{0,1}= \|(C_\mathrm{m}^{-1}\Delta_{\cS_\rho} - D\Delta_{\cS_R})u_1\|
\lesssim \delta\|u\|_{2,1}\,, \qquad (u\in H^{2,1})\,.
$$

To bound $\| (\lambda - A_R)^{-1} u \|_{2,1}$ we write $A_R = -(\lambda - A_R) + \lambda$ and apply the triangle inequality:
\begin{align*}
    \|(\lambda-A_R)^{-1}u\|_{2,1} &\lesssim \|A_R(\lambda-A_R)^{-1}u\|_{0,1} \\
    &\le \|u\|_{0,1} + |\lambda|\, \|(\lambda-A_R)^{-1}u\|_{0,1} \lesssim \|u\|_{0,1}\,,
\end{align*}
for all $u\in H^{0,1}$. In the last inequality we have used that $\|(\lambda-A_R)^{-1}\|_{0,1} \lesssim \min\{1, |\lambda|^{-1}\}$ by Lemma~\ref{lem:sectorial}.

To bound $\| (\lambda - A_{\rho})^{-1} u \|_{2,1}$, we use 
Eq.~\eqref{eq:domain-A-01}, Lemma~\ref{lem:sectorial}, and the fact that the resolvent $(\lambda - A_{\rho})^{-1}$ commutes with $A_{\rho}$, to see that
$$
\|(\lambda-A_\rho)^{-1}\|_{2,1}
\lesssim \|(\lambda-A_\rho)^{-1}\|_{2,1;\rho}
\lesssim \|(\lambda-A_\rho)^{-1}\|_{0,1;\rho}
\lesssim \min\{1, |\lambda|^{-1}\}\,.
$$

Applying the above estimates to Eq.~\eqref{eq:reso_id} we conclude that there exists a constant $C$ such that if $\delta:=R^{-1}\|\rho-R\|_{C^2}\le c$, where $c$ is the same constant as Lemma~\ref{lemma:difference}, then
\begin{equation}\label{eq:difference_resolvents}
    \| (\lambda - A_\rho)^{-1} - (\lambda - A_R)^{-1} \|_{2,1}
    \le C\delta \min\{1, |\lambda|^{-1}\}
\end{equation}
for all $\lambda$ with $\Re\,\lambda\ge -\frac{4-C_m r_{int}}{8} \nu\min\{1,\eps^{-\frac12}|\Im\, \lambda|\}$.

Therefore, applying Eq.~\eqref{eq:difference_resolvents} to Eq.~\eqref{eq:difference_semigroups} we arrive at
\begin{equation*}
    \|e^{tA_\rho}-e^{tA_R}\|_{2,1} \leq C \,\delta\int_0^\infty e^{-ts} \min\{1,s^{-1}\}\, ds\,.
\end{equation*}
For $t \geq 1$, the integral is uniformly bounded.
For $t< 1$, we have
$$
\int_0^\infty e^{-ts} \min\{1,s^{-1}\}\, ds
\le 1 + \int_1^{t^{-1}} \!\!\!s^{-1}\, ds +\int_{t^{-1}}^\infty 
\!t e^{-ts}\, ds
\le 2+\log t^{-1}\,,
$$
proving the claim.
\hfill $\Box$

\subsection{Nonlinear dynamics}\label{sec:nonlinear_dyn}

We extend the perturbation estimates from above to the nonlinear evolution generated by the FHN system on $\cS_{\rho}$. To prove Theorem~\ref{thm:persist} we need a perturbation estimate that controls the dependence of solutions to $\rho$. The size of the perturbation is measured in terms of the parameter $\delta:=R^{-1}\|\rho-R\|_{C^2}$.

\begin{proposition}\label{prop:continuous-rho} 
	Let $u\in C([0,T], H^{2,1})$ be a mild solution of Eq.~\eqref{eq:FHN} on $\cS = \cS_R$ with initial value $u\big\vert_{t=0} = u_0$. There are positive constants $\delta_*:=\min\{ c, \frac{1}{2C_0} \}$ and $C_0 := C' e^{\tilde{C} T}$ such that if $0<R\le 1$ and $\delta:= R^{-1}\|\rho-R\|_{C^2}\le \delta_*$, then the unique mild solution of Eq.~\eqref{eq:FHN} on $\cS_\rho$ with initial value $u_\rho\big\vert_{t=0} = u_0$ satisfies
	$$
	\sup_{0\le t\le T} \|u_\rho(t) - u(t)\|_{2,1}
	\le C_0\delta\,.
	$$
    The constant $c$ is as in Lemma~\ref{lemma:difference}, $C' = \sup_t \| e^{tA_{\rho}} \|_{2,1}$ and $\tilde{C} = C'C'_{\eta}$, with $C'_{\eta} = 3\eta^2 + 2\eta$.
\end{proposition}

\begin{proof}
	By Duhamel's formula
	\begin{equation}\label{eq:duhamel}
	u_\rho(t)-u(t) = (e^{tA_\rho}-e^{tA_R})u_0
	+ \int_0^t \left(
	e^{(t-s)A_\rho}N(u_\rho(s)) - e^{(t-s)A_R}N(u(s))\right)\, ds\,,
	\end{equation}
	so long as both solutions $u$ and $u_{\rho}$ exist. Using the triangle inequality and Eq.~\eqref{eq:semigroup_dif} we have
	\begin{align*} 
	\|e^{(t-s)A_\rho}N(u_\rho) &- e^{(t-s)A_R}N(u)\|_{2,1}
	\hskip -2.5cm\\
	&\le \|(e^{(t-s)A_\rho} - e^{(t-s)A_R})N(u)\|_{2,1} + 
	\|e^{(t-s)A_\rho}(N(u_\rho) - N(u))\|_{2,1} \\
	&\le C\delta (1 + \log (t^{-1} + 1))\|N(u)\|_{2,1} +
	C'\|N(u_\rho)-N(u)\|_{2,1}\,,
	\end{align*} 

    \noindent where $C'=\sup_t e^{tA_\rho}$. Next, we bound the nonlinearity. We have
        \begin{equation*}
            \| N(u) \|_{2,1} \leq C_{\eta} \| u \|_{2,1}
        \end{equation*}
        for all $u$ with $\| u_1 \|_{H^2} \leq \eta$, where $C_{\eta} := C_m^{-1}(\eta^2 + 4\eta)>0$. For the difference $N(u_{\rho}) - N(u)$ we have
        \begin{align*}
	           N_1(u_{\rho}) - N_1(u) &= -(u_{\rho_1}^3 - u_1^3) + (\alpha + 1) (u_{\rho_1}^2 - u_1^2) + (1 - C_\mathrm{m}^{-1})(u_{\rho_2} - u_2) \\
	           & =(u_{\rho_1} - u_1) \left( (\alpha + 1)(u_{\rho_1} + u_1) - (u_{\rho_1}^2 + u_{\rho_1}u_1 + u_1^2) \right) \\ 
	           & \qquad \qquad \qquad \qquad \qquad \qquad \quad + (1 - C_\mathrm{m}^{-1}) (u_{\rho_2} - u_2)\,,
	   \end{align*}    
       which implies 
       \begin{equation*}
           \| N(u_{\rho}) - N(u) \|_{2, 1} \leq C_\eta' \| u_{\rho} - u \|_{2,1}
       \end{equation*}
       for all $u_{\rho}$, $u$ with $\| u_{\rho_1} \|_{H^2}$, $\| u_1 \|_{H^2} \leq \eta$, where $C'_{\eta} := 3\eta^2 + 2\eta >0$.
       Hence,
       \begin{equation*}
           \|e^{(t-s)A_\rho}N(u_\rho) - e^{(t-s)A_R}N(u)\|_{2,1} \leq C C_\eta \delta (1+\log (t^{-1} + 1))\|u\|_{2,1} + C' C'_\eta \|u_\rho - u\|_{2,1}\,,
       \end{equation*}
       for positive constants $C$, $C'$, $C_{\eta}$ and $C'_{\eta}$ as defined above.

    Returning to Eq.~\eqref{eq:duhamel} we have    
	$$
	\|u_\rho(t)-u(t)\|_{2,1}
	\le C_T \delta \|u_0\|_{2,1} + \tilde{C}\int_0^t \|u_\rho(s)-u(s)\|_{2,1}\,ds\,,
	$$
	where $C_T = C(1+2 C_\eta T)(1+\log (T^{-1} + 1))$, and 
	$\tilde{C} := C' C'_\eta$.
	In the bound on the nonlinearity, we have used
	$\|u_\rho(t)\|_{2,1}\le \frac12 \eta$ for $0\le t\le T$. By Gr\"onwall's inequality
	$$
	\|u_\rho(t) - u(t)\|_{2,1} \le C_T \delta \| u_0 \|_{2,1} e^{\tilde{C} \int_0^tds} \le C_0\delta\|u_0\|_{2,1} \,,
	\qquad (0\le t\le T)
	$$
    where $C_0 := C_T e^{\tilde{C} T}$, and provided that $\sup_{0\le t\le T}\|u_\rho(t)\|_{2,1}\le \eta$. Since 
	$\|u(t)\|_{2,1}\le \frac12 \eta$ for $0\le t\le T$,
	by the triangle inequality 
	this is guaranteed by setting
	$\delta_*= \min\left\{c, \tfrac{\eta}{2C_0}\right\}$, where $c$ is the same as in Lemma~\ref{lemma:difference}. The statement of the proposition follows by absorbing $(1+2C_{\eta} T)(1+\log (T^{-1}+1))$ into $e^{\tilde{C}T}$.
	
\end{proof}

\subsubsection{Proof of Theorem~\ref{thm:persist}}
\label{subsec:proof-persist}

For the standard cylinder $\cS_R$, recall from Theorem~\ref{thm:nonli-stab} that there are constants $K_1 \geq 1$ and $\mu>0$ such that $\dist(u(t),\cM)\le K_1 e^{-\mu t} \|u_0-\Phi\|_{2,1}$ for every solution with initial values $u_0 \in \cU$. Here, $\cU$ is the neighbourhood of the traveling fronts $\Phi$ in $H^{2,1}$. By translation invariance
\begin{equation}
\label{eq:dist-M}
\dist(u(t),\cM)\le K_1 e^{-\mu t}\dist(u_0,\cM)
\end{equation}
for all $u$ with $u_0$ in some tubular neighborhood $\cW \subset \cM$ of the form $\cW = \left\{w\in H^{2,1}\ \big\vert \ \dist(w,\cM)<\xi\right\}$, for some $\xi>0$. By the triangle inequality and setting $T:= \mu^{-1}\log(2K_1)$ we have
$$
\sup_{0\le t\le T} \|u(t)\|
\le \sup_{0\le t\le T} \dist(u(t),\cM)+
\sup_{\Phi_h\in \cM} \|\Phi_h\|_{2,1} 
\le 2 \|\Phi\|_{2,1}
$$
for all solutions on the standard cylinder $\cS_R$ with $u_0 \in \cW$.

Turning to the warped cylinder $\cS_{\rho}$, let $u(t)$ be the mild solution of the reference system on $\cS_R$ with initial condition $u_0 \in \cW$. By Proposition~\ref{prop:continuous-rho}, there is a $\delta_0 > 0$ that is determined by $T$ and $C>0$ such that
\begin{equation}
\label{eq:delta}
\sup_{0\le t\le T}\|u_\rho(t)-u(t)\|_{2,1}\le C\delta \,,
\end{equation} 
provided that $\delta\le \delta_0$. Pick $\tilde{\delta}:=\min\{\delta_0,\frac{\eta}{2C}\}$. From Eq.~\eqref{eq:dist-M} and Eq.~\eqref{eq:delta}, and assuming that $\delta \leq \tilde{\delta}$, we have
\begin{equation}\label{eq:proof-persist-interval}
    \begin{split}
        \dist(u_\rho(t),\cM) &\le \dist(u(t),\cM) + \|u_\rho(t)-u(t)\|_{2,1}\\
        & \le K_1 e^{-\mu t} \dist(u_0,\cM) + C\delta\\
        &< 2K_1\eta
    \end{split}
\end{equation}
for all $t\in[0, T]$. Equation~\eqref{eq:proof-persist-interval} implies that $\|u_\rho(t)\|_{2,1}\le 2\|\Phi\|_{2,1}$. Moreover, $u_{\rho}(T) \in \cW$ whenever $u_{\rho}(0)\in \cW$ because
$$
\dist(u_\rho(T),\cM) \le \frac12 \dist(u_0,\cM) + C\delta< \nu\,.
$$
By induction
\begin{equation*}
    \dist(u_{\rho}((k+1)T), \cM)) \leq \frac{1}{2} dist (u_{\rho}(kT), \cM) + C\delta\,,\qquad (k\in \N)\,.
\end{equation*}
Solving the recursion, we conclude that $\dist(u_{\rho}(kT), \cM) \leq 2^{-k} \dist(u_0, \cM) + \left( \sum_{j=0}^{k-1} 2^{-j} \right) K_1 \delta$ and by Eq.~\eqref{eq:proof-persist-interval}
\begin{equation*}
    \dist(u_\rho(t),\cM) \le 2^{-k} C_0\, \dist(u_0,\cM) + (2+C_0)C\delta
\end{equation*}
for all $t$ with $kT \leq t \leq (k+1)T$ and all $k \in \mathbb{N}_0$.

Using that $T = \mu^{-1} ln(2 K_1)$ we conclude that
\begin{equation*}
    \dist(u_{\rho}(t), \cM) \leq K'_1 e^{-\nu t} \dist(u_0, \cM) + K'_2 \delta
\end{equation*}
where $K'_1 := 2K_1$, $K'_2 := (2 + K_1) C_0$, $\nu = \frac{\mu ln2}{ln(2 K_1)}$, $K_1$ is from Theorem \ref{thm:nonli-stab} and $C_0 := C' e^{\tilde{C}T}$, as in Proposition \ref{prop:continuous-rho}.

\hfill $\Box$

\section{Numerical simulations}\label{sec:numerics}

In this section, we present the numerical analysis of Eq.~\eqref{eq:FHN} on $S = S_{\rho}$. The equation was solved using the method of lines: space was discretized with finite differences, and the resulting stiff system of ordinary differential equations was integrated in time using MATLAB’s implicit solver.

We considered a rectangular domain in $(x,\theta)$, with $x \in [0,L]$ and $\theta \in [0,2\pi]$. Neumann (zero-flux) boundary conditions were imposed in the $x$-direction, while periodic boundary conditions were applied in the $\theta$-direction. The domain was discretized using a uniform grid. The corresponding mesh points are given by
\begin{equation*}
    x_i= i \Delta x, \;\;\; i=0, \ldots, N-1, \;\;\; \text{and} \;\;\;
    \theta_j = j \Delta \theta, \;\;\; j=0, \ldots, M-1 \,,
\end{equation*}
where $\Delta x = \frac{L}{N}$ and $\Delta \theta = \frac{2\pi}{M}$.

\paragraph{Numerical scheme of axial direction.} The operator in the longitudinal coordinate $x$ has the divergence form
\begin{equation*}
    \frac{1}{w_1(x)} \partial_x\left(w_2(x) \partial_x(x)\right)    
\end{equation*}
where
\begin{equation*}
    w_1(x) = \rho(x) \sqrt{1+\rho'(x)^2} \;\;\; \text{and} \;\;\; w_2(x) = \frac{\rho^3(x)}{\sqrt{1+\rho'(x)^2}}\,.
\end{equation*}
The discrete approximation is constructed by composing two centered finite difference operators. First, the derivative of $u$ is approximated by a second-order central difference,
\begin{equation*}
    \left( D_x^1 u \right)_i = \frac{u_{i+1} - u_{i-1}}{2 \Delta x}, \;\;\; 2 \leq i \leq N-1,
\end{equation*}
with one-sided second-order stencils applied at the boundaries to enforce the Neumann condition. A second central-difference operator is then applied to $q_i = w_{2,i}(D_x^1 u)$:
\begin{equation*}
    (D_x^2 q)_i = \frac{q_{i+1} - q_{i-1}}{2 \Delta x}\,,
\end{equation*}
which is the approximation of $\partial_x(w_2 \partial_x u)$. Expanding the stencil for interior points gives
\begin{equation*}
    \frac{w_{2,i+1}}{4(\Delta x)^2 w_{1,i}} u_{i+2} - \frac{w_{2,i+1} + w_{2, i-1}}{4(\Delta x)^2 w_{1,i}} u_{i} + \frac{w_{2,i-1}}{4(\Delta x)^2 w_{1,i}} u_{i-2}\,.
\end{equation*}

\paragraph{Numerical scheme of angular direction.} The standard second-order central difference stencil was used
\begin{equation*}
    \partial_{\theta}^2 u \big\vert_{\theta_j} \approx \frac{u_{j-1} - 2u_j + u_{j+1}}{(\Delta \theta)^2}, \;\;\; j=1, \ldots, M
\end{equation*}
with periodic boundary conditions.

\vspace{0.25cm}

Equation~\eqref{eq:FHN} was numerically solved on the standard cylinder $S_{R}$ using the same uniform meshes as described above. Spatial derivatives in both the axial and angular directions were approximated with standard second-order central finite differences.

\vspace{0.25cm}

\begin{figure}[tbhp]%[ht]%[tbhp]
	\centering
	\includegraphics[width=1\linewidth]{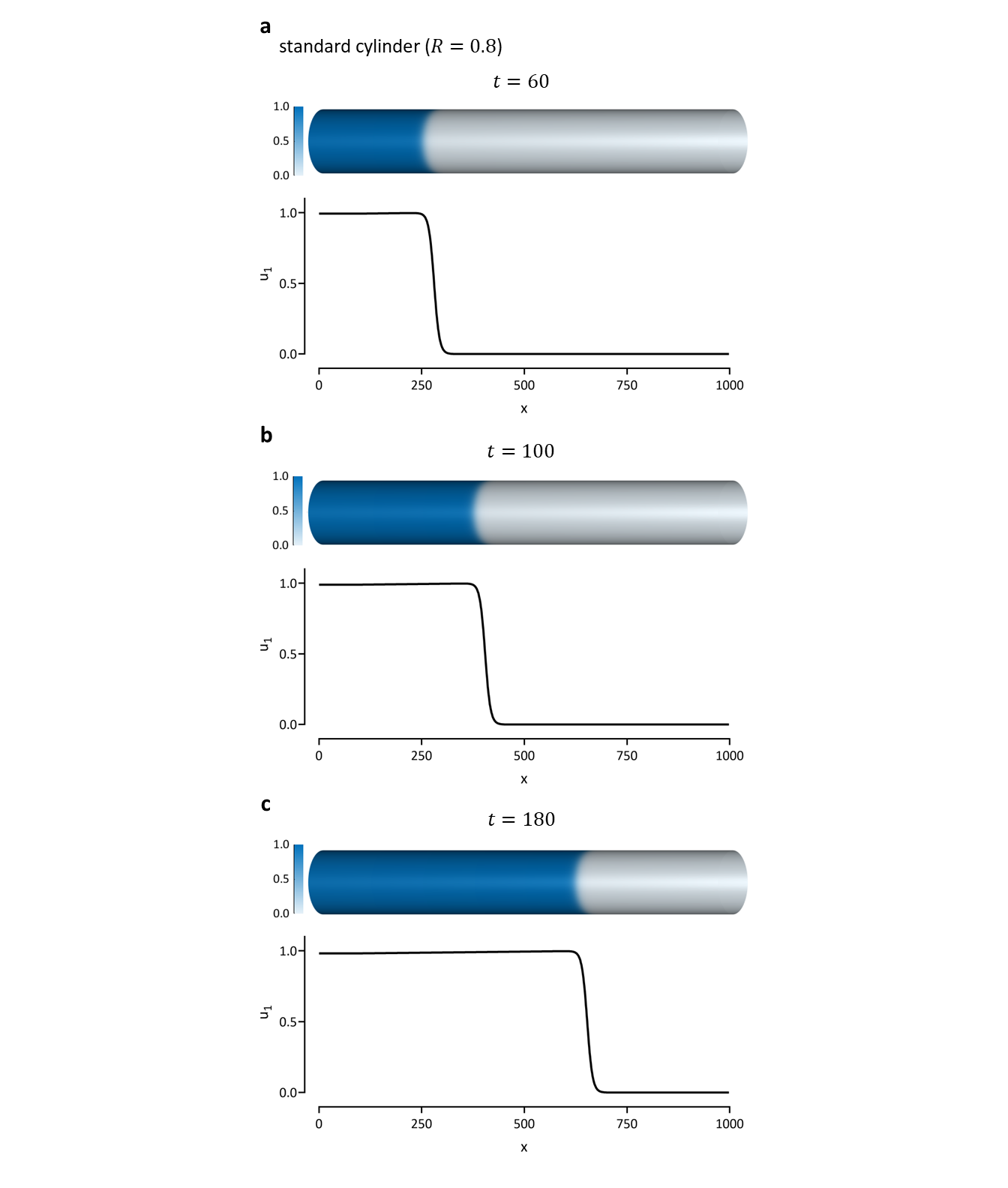}
	\caption{\textbf{Traveling fronts on the surface of a constant cylinder.} \textbf{a-c.} \textit{Top}: Heat maps on three-dimensional cylindrical surfaces illustrating the geometries considered. \textit{Bottom:} Snapshots of traveling front solutions at three time instances, $t = 60, 100, 180$. Note that the front preserves its shape as it propagates along the cylinder. The parameter values are $\alpha = 0.01$, $\eps = 0.0001$, $\gamma = 7$, $C_m = 1$, $r_{\mathrm{int}} = 0.1$, and $L = 1000$.
	} 
	\label{fig:constant}
	
\end{figure}

Figure~\ref{fig:constant}a-c illustrates the propagation of a traveling front on a cylindrical surface $\cS_R$ with fixed radius $R=0.8$. Each panel shows a three-dimensional representation of the cylinder together with the corresponding profile of the excitation variable $u_1$ as a function of the spatial coordinate $x$. The snapshots, taken at successive times, reveal the evolution of the front. The results demonstrate that the front remains stable and propagates without altering its shape.

\begin{figure}[tbhp]%[ht]%[tbhp]
	\centering
	\includegraphics[width=1\linewidth]{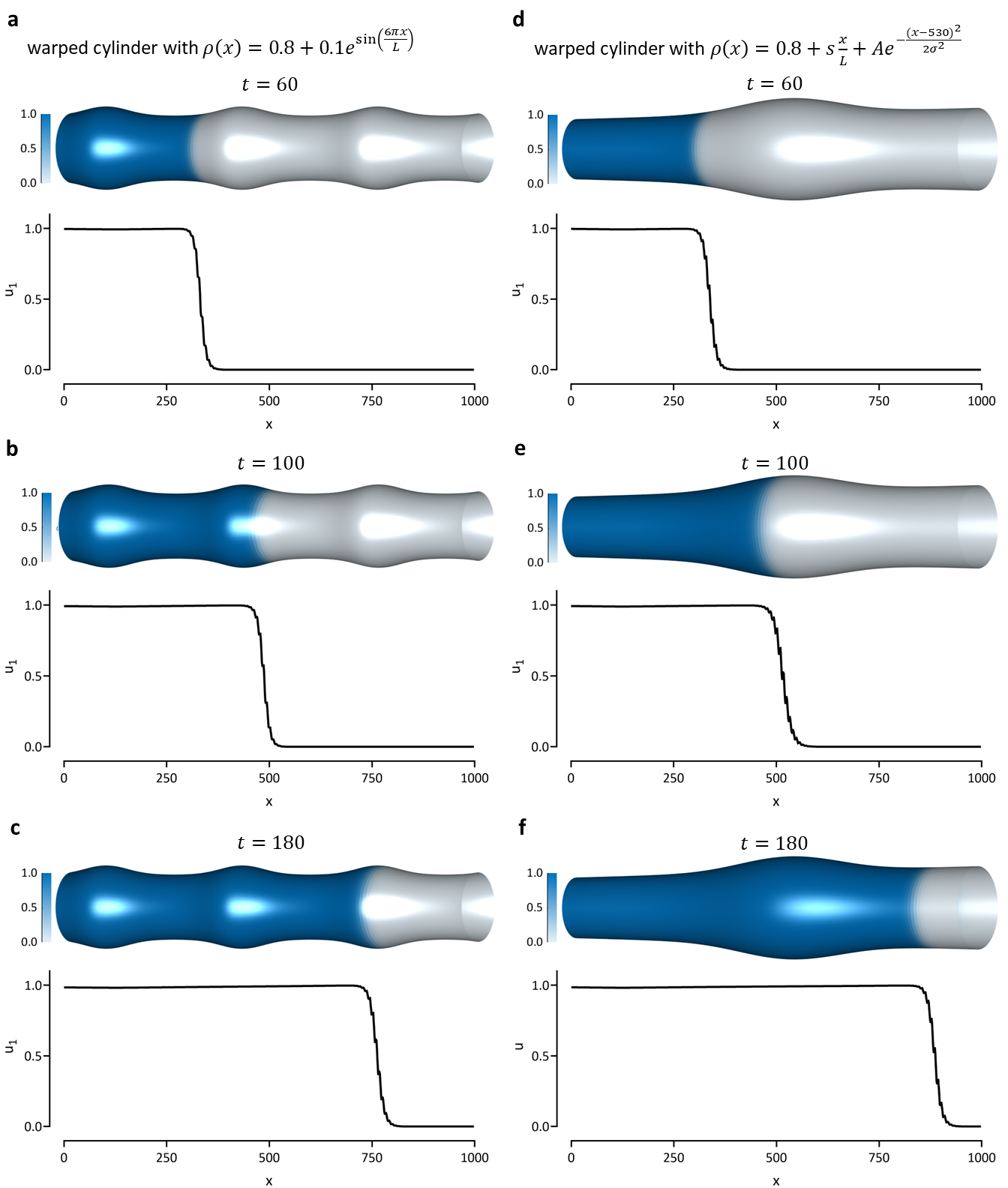}
	\caption{\textbf{Traveling fronts on warped cylinders.} \textbf{a-c.}The cylinder radius is given by $\rho(x) = 0.8 + 0.1 e^{sin(\frac{6 \pi x}{L})}$, illustrating an example of a pearl-on-a-string morphology of neuronal axons. \textbf{d-f.} The cylinder radius is given by $\rho(x) = 0.8 + s\frac{x}{L} + A e^{-\frac{(x-x_0)^2}{2\sigma^2}}$, representing a localized axonal swelling. The parameter values are the same as in Fig.~\ref{fig:constant}. The near-front solutions remain stable but propagate with different speeds, highlighting the influence of axonal morphology on signal transmission.
	} 
	\label{fig:warped}
	
\end{figure}

As a first example of a warped cylindrical surface $\cS_{\rho}$, we consider the so-called pearls-on-a-string axonal morphology. The radius is defined as
\begin{equation}\label{eq:pearl}
    \rho(x) = 0.8 + 0.1 e^{sin(\frac{6 \pi x}{L})}\,,
\end{equation}
where $L$ denotes the total length of the axon. The pearls-on-a-string morphology refers to axonal structures in which segments of narrow cable alternate with periodic bulges or “pearls” \cite{pearl}. The size and number of these pearls can vary, enabling axons to regulate the velocity of electrical signal propagation. Figure~\ref{fig:warped}a-c illustrates the profile of the excitation variable $u_1$ along the spatial coordinate $x$, together with three-dimensional renderings of the warped cylindrical surfaces whose radius is defined by Eq.~\eqref{eq:pearl}. The snapshots correspond to the same time instants $t$
as those shown in Fig.~\ref{fig:constant}. For these times, the position of the traveling front along $x$ differs between the two cylinders, indicating that the propagation speed is influenced by the surface geometry.

Another example of a warped cylindrical surface is given by axons exhibiting a single swelling -- a focal enlargement that can occur in both healthy and diseased axons \cite{Kolaric}. We model such an axon by introducing a Gaussian bump in the radius profile. In this case, the radius is defined as
\begin{equation}\label{eq:swelling}
    \rho(x) = 0.8 + s\frac{x}{L} + A e^{-\frac{(x-x_0)^2}{2\sigma^2}}\,,
\end{equation}
where the second term represents a small linear slope, and the third term corresponds to a Gaussian bump centered at $x_0$. The parameter values used in the numerical simulations are: $s=0.3$, $L=1000$, $A=0.4$, $x_0 = 530$, and $\sigma=120$. The corresponding profile of $u_1$, and the axonal representations of $\cS_{\rho}$ with $\rho$ defined by Eq.~\eqref{eq:swelling}, are shown in Fig.~\ref{fig:warped}d-f. The traveling wave remains stable and preserves its front-like shape as it propagates along the axon.

\section{Discussion}\label{sec:discussion}

In this work, we analyzed the stability of traveling front solutions of the FitzHugh – Nagumo system on thin, infinitely long cylindrical surfaces. For cylinders of constant radius, the results show that traveling front solutions are nonlinearly stable under perturbations on the initial condition. When the cylinder’s radius varies slowly along its axis, thereby modeling more realistic axonal geometries, the fronts adapt to the local geometry, and solutions that begin near a traveling front remain close to the corresponding family of propagating fronts for all time. The first result is established using an exponentially decaying estimate for the semigroup generated by the linearized operator of the FitzHugh-Nagumo system in the moving frame. The second result is obtained by viewing the dynamics on a warped cylinder as a perturbation of the flow on a standard cylinder. The comparison of the corresponding linearized flows is derived from semigroup estimates that rely on the analyticity and uniform boundedness of the semigroups. Gronwall’s inequality is then employed to extend these estimates to the nonlinear flow. Combined with the exponential decay of perturbations established for the standard cylinder, this argument yields the persistence of traveling fronts on warped cylinders. The drawback of this method is that Gronwall’s inequality produces constants with unfavorable $\eps$-dependence.

Our numerical simulations were motivated by recent developments in theoretical neuroscience. A new study \cite{pearl} suggests that pearls-on-a-string axonal morphologies are a common feature of healthy unmyelinated axons, where they play a regulatory role in modulating the propagation speed of action potentials. This speed is, in turn, a key determinant of effective neural communication. Various factors -- such as the size and spacing of the pearls, as well as the permeability of the axonal membrane -- may significantly influence signal propagation. We aim to investigate the impact of axonal geometry on the dynamics and speed of signal propagation in subsequent work.

We also performed numerical simulations on an axon with a single swelling. Axonal swellings are abnormal enlargements that can develop along the length of an axon \cite{Kolaric}. Depending on their size, such swellings are often pathological and constitute early indicators of neuronal dysfunction and axonal degeneration \cite{swelling_Johnson, swelling_Maia}. The formation of swellings not only alters the speed of signal propagation but can also lead to partial or complete blockage of electrical signals. In future work, we aim to determine the conditions under which these blockages occur.

Other possible extensions of this work include the following. (i) The stability analysis of periodic traveling waves on cylindrical surfaces similar to those considered above. Such waves typically take the form of pulse trains, arising from the continuous firing of the neuron. (ii) The study of the existence and stability of traveling pulses along axons with heterogeneous myelination. Myelin is a fatty sheath that surrounds axons and increases the speed of electrical signal propagation. In general, myelinated axons maintain a more regular cylindrical shape than unmyelinated ones, however, swellings may also occur, particularly as a consequence of neurodegenerative diseases or trauma \cite{swelling_Johnson}.

%\newpage

\paragraph{Acknowledgments.} The author would like to thank Almut Burchard and Georgia Karali for helpful comments on the manuscript, and Israel Michael Sigal and Wilten Nicola for stimulating discussions. This work was supported in part by the Fields Institute for Research in Mathematical Sciences.

%%%%%%%%%%%%%%%%%%%%%%%%%
%%%%%%%%%%%%%%%%%%%%%%%%%

\newpage

\appendix

\section{Domain characterization of \( L \)}\label{app:equivalence}

\begin{lemma}
	\label{lem:L-domain}
	Let \( L \) be the linear operator defined in Eq.~\eqref{eq:L}. Then the domain of \( L \) coincides with the function space \( H^{2,1} \), and the associated graph norm satisfies the equivalence:
	\[
	\|u\|_{2,1} \lesssim \|Lu\| + \|u\| \lesssim \|u\|_{2,1}, \qquad \text{for all } u \in H^{2,1}.
	\]
\end{lemma}

\begin{proof}
	We begin by estimating the longitudinal derivative of \( u_1 \) using a standard interpolation inequality:
	\[
	\|\partial_z u_1\| \leq \frac{1}{2c} \|\Delta_{\mathcal{S}_R} u_1\| + \frac{c}{2} \|u_1\|.
	\]
	Next, define
	\begin{equation}
	\label{eq:b}
	b := \sup_{z \in \mathbb{R}} \left| f'(\phi_1(z)) - f'(0) \right| < \infty,
	\end{equation}
	noting that the nonlinearity \( f \) is smooth and \( \phi_1 \) is bounded.
	Using the reverse triangle inequality and applying the definition of \( L \), we obtain a lower bound
	\[
	\|Lu\| \geq \min\left\{ \frac{1}{2},\, c \right\} \left\| \left( \Delta_{\mathcal{S}_R} u_1, \partial_z u_2 \right) \right\| - \left( \frac{c^2}{2} + b \right) \|u_1\|.
	\]
	Hence, we derive the estimate
	\begin{align*}
	\left(1 + \frac{c^2}{2} + b \right)( \|Lu\| + \|u\| )
	&\geq \|Lu\| + \left(1 + \frac{c^2}{2} + b \right) \|u\| \\
	&\geq \min\left\{ \frac{1}{2},\, c \right\} \|u\|_{2,1}.
	\end{align*}
	For the upper bound, a direct application of the triangle inequality gives
	\[
	\|Lu\| + \|u\| \leq \max\left\{ \frac{3}{2},\, c,\, 1 + \frac{c^2}{2} + b \right\} \|u\|_{2,1},
	\]
	as required.
\end{proof}

\section{Dissipativity and generation}\label{app:dissipativity}

\begin{lemma}
	\label{lem:dissipative}
	Let \( B \) be a closed, densely defined dissipative operator on a Hilbert space \( H \). Then \( B \) generates a strongly continuous semigroup of contractions \( e^{tB} \). Moreover, the spectrum of \( B \) lies entirely within the left half-plane \( \{ \lambda \in \mathbb{C} \mid \Re \lambda \leq 0 \} \), and the resolvent satisfies the estimate
	\[
	\| (\lambda - B)^{-1} \| \leq \frac{1}{\Re \lambda}, \qquad \text{for all } \Re \lambda > 0.
	\]
\end{lemma}

\begin{proof}
	By dissipativity, the operator \( I - B \) is injective, and for all \( v \in D(B) \):
	\[
	\| (I - B)v \| \geq \|v\|^{-1} \left| \Re \langle (I - B)v, v \rangle \right| \geq \|v\|.
	\]
	A similar argument for the adjoint $B^*$ implies that $\Ran( I - B )$ is dense in \( H \).
	
	Let \( w_0 \in H \) be arbitrary. There exists a sequence \( (w_n) \subset \operatorname{Ran}(I - B) \) with \( w_n \to w_0 \). For each \( n \), let \( v_n \in D(B) \) satisfy \( (I- B)v_n = w_n \). Then,
	\[
	\|v_n - v_m\| \leq \|(I - B)(v_n - v_m)\| = \|w_n - w_m\|,
	\]
	so \( (v_n) \) is a Cauchy sequence and converges to some \( v_0 \in H \). Since \( B \) is closed, \( w_0 = (I - B)v_0 \), showing that \( I - B \) is surjective.
	
	By the Lumer–Phillips theorem, \( B \) generates a strongly continuous semigroup of contractions. The resolvent bound then follows from the Hille–Yosida theorem \cite[Theorem 1.3.1]{Pazy}.
\end{proof}

\section{Stability theorem}\label{app:stability_thm}

\begin{theorem}[Theorem 4.3.5, \cite{KPbook}]
    Consider the nonlinear problem
    \begin{equation*}
        \partial_t u= \cF(u)\,,
    \end{equation*}
    where $\cF:Z \subset X \mapsto X$ has an $N$-fold symmetry. Suppose that $\cF(\Phi) = 0$, and that the linearization, $L=\nabla_u \cF(\phi)$, and nonlinearity, $\mathcal{N}$, defined as
    \begin{equation*}
        \cN (v) := \cF(\phi + v) - Lv,
    \end{equation*}
    satisfy the following hypothesis, for some $\sigma>0$:
    \begin{enumerate}[(a)]
        \item There is a constant $\sigma > 0$ such that
        \begin{equation}
            S_+ := \sigma(L)\cap \left\{ \lambda \in \mathbb{C}: Re \lambda \geq -\sigma \right\} = \left\{ 0 \right\}\,.
        \end{equation}
        %%%%%
        \item The generalized kernel of $L$ is spanned by the symmetry eigenfunctions, with
        \begin{equation}
            \text{gker}(L) = \text{ker}(L) = \text{span}\left\{T_1'\phi, \ldots, T'_N \phi \right\}\,.
        \end{equation}
        %%%%%
        \item There exists $M>0$ such that the spectral projection $\Pi_+:X \mapsto \text{ker}(L^a)^\perp$ complimentary to the set $S_+$ satisfies the resolvent estimate
        \begin{equation}
            \| R(\lambda; L)\Pi_+ f \|_{Y} \leq M \| f \|_Y, \;\;\; Re \lambda \geq -\sigma \,.    
        \end{equation}
        %%%%%
        \item The gradient $\nabla_u\cF$ is locally Lipschitz on bounded sets; that is, for each $R>0$ there exists $M>0$ such that
        \begin{equation}
            \left\| (\nabla_u \cF(u) - \nabla_u \cF(v))w \right\|_Y \leq M \| u - v \|_Y\, \| w \|_Y\,,
        \end{equation}
        as long as $\| u \|_Y$, $\| v \|_Y \leq R$.
        %%%%%
        \item The nonlinearity $\cN$ is quadratic in $\| \cdot \|_Y$ near zero, i.e., there exists $R$, $M>0$ such that
        \begin{equation}
            \| \cN(v) \|_Y \leq M \| v \|_Y^2\,,
        \end{equation}
        for all $\| v \|_Y \leq R.$
    \end{enumerate}
    Then for any $\tilde{\sigma} \in (0, \sigma)$ the manifold $\cM_T$ of equilibria is asymptotically orbitally stable in $\| \cdot \|_Y$ with exponential rate $\tilde{\sigma}$.
    
\end{theorem}

\end{document}